\documentclass{svmult}
\usepackage{makeidx}         
\usepackage{graphicx}        
\usepackage{multicol}        
\usepackage[bottom]{footmisc}

\makeindex             

\newtheorem{thm}{Theorem}[section]
\newtheorem{lem}{Lemma}[section]
\newtheorem{prop}{Proposition}[section]
\newtheorem{cor}{Corollary}[section]
\newtheorem{rmk}{Remark}[section]
\newtheorem{exap}{Example}[section]
\newtheorem{defi}{Definition}[section]

\begin{document}

\title*{Two-parameter $p, q$-variation Paths and Integrations of Local Times }
\titlerunning{Two-parameter $p, q$-variation paths and integrations of local times }
\author{Chunrong Feng\inst{1,2},
Huaizhong Zhao\inst{1}}
\authorrunning{C. Feng and H. Zhao}
\institute{ Department of Mathematical Sciences, Loughborough
University, LE11 3TU, UK. \texttt{C.Feng@lboro.ac.uk},
\texttt{H.Zhao@lboro.ac.uk}
\and School of Mathematics and System Sciences, Shandong University, Jinan, Shandong Province, 250100, China}
\maketitle
\newcounter{bean}
\begin{abstract}
In this paper, we prove two main results. The first one is to give
a new condition for the existence of two-parameter $p,
q$-variation path integrals. Our condition of locally bounded
$p,q$-variation is more natural and easy to verify than those of
Young. This result can be easily generalized to multi-parameter
case. The second result is to define the integral of local time
$\int_{-\infty}^\infty\int_0^t g(s,x)d_{s,x}L_s(x)$ pathwise and
then give generalized It$\hat {\rm o}$'s formula when
$\nabla^-f(s,x)$ is only of bounded $p,q$-variation in $(s,x)$. In
the case that $g(s,x)=\nabla^-f(s,x)$ is of locally bounded
variation in $(s,x)$, the integral $\int_{-\infty}^\infty\int_0^t
\nabla^-f(s,x)d_{s,x}L_s(x)$ is the Lebesgue-Stieltjes integral
and was used in Elworthy, Truman and Zhao \cite{Zhao}. When
$g(s,x)=\nabla^-f(s,x)$ is of only locally $p, q$-variation, where
$p\geq 1$,$q\geq 1$, and $2q+1>2pq$, the integral is a
two-parameter Young integral of $p,q$-variation rather than a
Lebesgue-Stieltjes integral. In the special case that
$f(s,x)=f(x)$ is independent of $s$, we give a new condition for
Meyer's formula and $\int_{-\infty}^\infty L_t(x)d_x\nabla^-f(x)$
is defined pathwise as a Young integral. For this
we prove the local time $L_t(x)$ is of $p$-variation in $x$ for
each $t\geq 0$, for each $p>2$ almost surely ($p$-variation in the
sense of Lyons and Young, i.e. $\sup\limits_{E: \ a\ finite\  partition\ of\
[-N,N]} \sum\limits_{i=1}^m|L_t(x_i)-L_t(x_{i-1})|^p<\infty$).

 \vskip 5pt
\end{abstract}
 \renewcommand{\theequation}{\arabic{section}.\arabic{equation}}
 \section{Introduction}

The classical It$\hat {\rm o}$'s formula for twice differentiable functions has played a central role in stochastic analysis and almost all aspects of its applications and connection with analysis, PDEs, geometry, dynamical systems, finance and physics. But the restriction of It$\hat {\rm o}$'s formula to functions with twice differentiability often encounter difficulties in applications. Extensions to less smooth functions are useful in studying many problems such as partial differential equtions with some singularities and mathematics of finance. Generally speaking, for any absolutely continuous function whose derivative $f'$ exists almost everywhere, and a continuous semi-martingale $X_t$, there exists $A_t$ such that
 \begin{eqnarray}
 f(X_t)=f(X_0)+\int_0^t f'(X_s)dX_s+A_t
 \end{eqnarray}
 and for the time dependent case, the corresponding formula is
 \begin{eqnarray}
f(t,X_t)=f(0,X_0)+\int_0^t {\partial\over \partial s}f(s,X_s)ds+\int_0^t\nabla f(s,X_s)dX_s+A_t.
 \end{eqnarray}
 To find $A_t$ in both cases especially a pathwise formula becomes key to establish a useful extension to It$\hat {\rm o}$'s formula. In fact investigations already began in Tanaka \cite{tan} with a beautiful use of local times introduced in L\'evy \cite{levy}. The generalized It$\hat {\rm o}$'s formula in one-dimension for time independent convex functions was developed in Meyer \cite{meyer} and for superharmonic functions in multidimensions in Brosamler \cite{brosamler} and for distance function in Kendall \cite{kendall} and more recently for time dependent functions in Peskir \cite{peskir3}, Ghomrasni
 and Peskir \cite{Peskir1}  and Elworthy, Truman and Zhao \cite{Zhao}. Meyer \cite{meyer} proved if $f$ is a convex function (or difference of two convex functions), then
\begin{eqnarray}\label{zhao2}
f(X_t)=f(X_0)+\int_0^t \nabla^- f(X_s)dX_s+\int_{-\infty}^\infty L_t(x)d_x\nabla^- f(x) \ a.s.
\end{eqnarray}
where $\nabla ^-f(x)$ is of bounded variation and
$\int_{-\infty}^\infty L_t(x)d_x\nabla ^-f(x)$ is a
Lebesgue-Stieltjes integral associated with the measure $d_x\nabla
^-f(x)$. Elworthy, Truman and Zhao \cite{Zhao} proved if
$f(t,x)=f_h(t,x)+f_v(t,x)$, where $\Delta^- f_h(t,x)$ and
$\nabla^- f(t,x)$ exist and are left continuous, and
$\nabla^-f_v(t,x)$ is of locally bounded variation in $x$ for a
fixed $t$ and of locally bounded variation in $(t,x)$, then
 \begin{eqnarray}\label{zhao1}
&&f(t,X(t))-f(0,X(0))\nonumber\\
&=&\int _0^t{\partial ^-\over \partial s} f(s,X(s)){\rm d}s+\int
_0^t\nabla^- f(s,X(s))dX_s\nonumber\\
&&+{1\over 2} \int_0^t \Delta^- f_h(s,X(s))d<\hskip-4pt X\hskip-4pt>_s
+
\int _{-\infty}^{\infty }L _t(x){\rm d}_x\nabla^- f_v(t,x)\nonumber\\
&&-\int _{-\infty}^{+\infty}\int _0^{t}L _s(x){\bf \rm d}_{s,x}\nabla^- f_v(s,x) \ \ a.s.
\end{eqnarray}
where $\int _{-\infty}^{+\infty}\int _0^{t}L _s(x){\bf \rm
d}_{s,x}\nabla^- f_v(s,x)$ is a space-time Lebesgue-Stieltjes
integral and needless to say, defined pathwise.
Elworthy-Truman-Zhao's formula was given in a very general form.
It includes as special cases classical It$\hat {\rm o}$'s formula,
Tanaka's formula, Meyer's formula, Azema-Jeulin-Knight-Yor's
formula \cite{yor1}. A special and earlier version of
Elworthy-Truman-Zhao's formula was obtained by Peskir
\cite{peskir2} independently. Feng and Zhao \cite{feng} extended
(\ref{zhao1}) to two dimensions. Noticing that the nonexistence of
local time in two dimensions gives an essential difficulty in
extending (\ref{zhao1}) to 2-dimensions, so the extension was
nontrivial and the key was to define the stochastic
Lebesgue-Stieltjes integral.

On the other hand, there are some works to define
 $\int_{-\infty}^\infty \nabla^- f(x) d_xL_t(x)$ and
 $\int_{-\infty}^\infty\int_0^t \nabla^- f(s,x)d_{s,x}L_s(x)$ in $L^2(dP)$ or
 in $L^1(dP)$  in connection with It$\hat {\rm o}$'s formula by using terms in
 (\ref{zhao2}) other than the last term or backward-forward stochastic integrals
 (Bouleau and Yor \cite{bou}, Eisenbaum \cite{eisenbaum1},
\cite{eisenbaum2}, Flandoli, Russo, Wolf \cite{frw}, F\"ollmer and
Protter \cite{Protter2}, Moret and Nualart \cite{nualart}, etc.)
and the work of Rogers and Walsh \cite{rog} using excursion
fields. Generally speaking, one expects stronger conditions for
the pathwise existence of the integrals of local times. However,
in the framework of classical integrals with respect to measures, locally bounded variation
in $x$ for fixed $t$ and locally bounded variation in $(t,x)$ are
minimal conditions on $\nabla^-f(t,x)$ to generate a measure, so
it seems impossible to go beyond Elworthy-Truman-Zhao's formula.
We remark that the striking fact that $L_t(x)$ is of bounded
quadratic variation in $x$ in the sense of Revuz and Yor
\cite{yor} and increasing in $t$ did not play a significant role
in the proof of (\ref {zhao1}). It is therefore reasonable to
conjecture that the conditions in \cite{Zhao} defining the
integrals of local times pathwise can be weakened. Inevitably, we
have to go beyond Lebesgue-Stieltjes integral as it seems to us that
Elworthy-Truman-Zhao's formula has achieved the best in the
Lebesgue-Stieltjes integral framework. Here we use Young's idea of
integration (Lyons \cite{terry2}, \cite{terry1}, Lyons and Qian
\cite{terry}, Young \cite{young1}, \cite{young}) to define the
integral of local time to go beyond the bounded variation
condition. We would like to remark that the quadratic variation in
the sense of Revuz and Yor is not enough to define Young integral
for local times. So it is necessary to prove local time $L_t(x)$ is of bounded
$p$-variation for each $p>2$ in the sense of Young almost surely.
The main difficulty is overcome by using the idea of controlling
the $p$-variation of continuous paths via the variations through
dyadic partitions. This idea was originated by L\'evy, used in
\cite{bass}, \cite{hambly}, \cite{ledoux} to prove the  Brownian path is
of bounded $p$-variation for $p>2$.

Using Young's integration of one parameter $p$-variation, we can
immediately define $\int_{-\infty}^\infty \nabla^-f(x) d_xL_t(x)$
as a Young integral  if $\nabla^-f(x)$ is of bounded
$q$-variation ($1\leq q<2$). Then a new extension of
Meyer's formula to $f$ where $\nabla^-f(x)$ is of bounded
$q$-variation ($1\leq q<2$) follows immediately. However one
can immediately realize the difficulty of defining the
two-parameter integral when we work on time dependent $f$. Young
\cite{young} considered this problem, but his conditions are
strong and difficult to check. It seems to us that the theory of
two-parameter $\Phi_1, \Psi_1$-variation ($p,q$-variation as a
special case) integration has not been investigated and developed
well in the literature. Inspired by the work of Young \cite{young} and
Lyons and Qian \cite{terry}, in this paper,  we give a new
condition for the existence of two-parameter Young integral
(Theorem \ref{ffcr1}). We consider a continuous function $F(x,y)$
being of bounded $\Phi$-variation in $x$ uniformly in $y$, and
being of bounded $\Psi$- variation in $y$ uniformly in $x$;
$G(x,y)$ being of bounded $\Phi_1, \Psi_1$-variation in ($x,y$),
i.e.
\begin{eqnarray}\label{cond1}
\sup_{E\times E'}\sum_{j=1}^{N'}\Psi_1\left(\sum_{i=1}^N\Phi_1(|\Delta_i\Delta_j G|)\right)<\infty,
\end{eqnarray}
where $E\times E':=\{x'=x_0<x_1<\cdots<x_N=x'', y'=y_0<y_1<\cdots<y_{N'}=y''\}$ is an arbitrary partition of $[x',x'']\times[y',y'']$, $|\Delta_i\Delta_j G|=|G(x_i,y_j)-G(x_{i-1},y_j)-G(x_i,y_{j-1})+G(x_{i-1},y_{j-1})|$, and
 $\Psi_1$, $\Phi_1$ are convex functions. Then if there exist monotone increasing functions $\varrho$ and $\sigma$ subject to $\varrho(u)\sigma(u)=u$ such that
\begin{eqnarray}
\sum\limits_{m,n}\varrho[\varphi({1\over n})]\sigma[\psi({1\over m})]\varphi_1[{1\over n}\psi_1({1\over m})]<\infty,
\end{eqnarray}
 then the integral
\begin{eqnarray}
&&\int_{x'}^{x''}\int_{y'}^{y''}F(x,y)d_{x,y}G(x,y)\nonumber\\
&=&\lim_{m(E\times E')\to 0}\sum_{i=1}^N\sum_{j=1}^{N'}F(x_{i-1},y_{j-1})\Big(G(x_i,y_j)-G(x_{i-1},y_j)\nonumber\\
&& \hskip 3cm-G(x_i,y_{j-1})+G(x_{i-1},y_{j-1})\Big)
\end{eqnarray}
is well defined. Here $\phi, \psi, \phi_1, \psi_1$ are inverse
functions of $\Phi, \Psi, \Phi_1, \Psi_1$, respectively. For this
we use Lyons' idea of control function to two-parameter case. We
also prove a dominated convergence theorem (Theorem \ref{our3})
for the integral. Then we apply this to establish the integral of
local time $\int_{-\infty}^\infty \int_0^t L_s(x) d_{s,x} \nabla
^-f(s,x)$ pathwise, where $\nabla ^-f(t,x)$ is of locally bounded
p,q-variation with $p,q \geq 1$, $2q+1>2pq$. This is new in
the literature. Under this condition we establish generalized It$\hat
{\rm o}$'s formula with the help of the dominated convergence
theorem. We believe our results of the two-parameter
$p,q$-variation path integration are new and has independent
interest.

To compare our condition (\ref{cond1}) with that of Young, we
quote his condition here
\begin{eqnarray}\label{cond2}
|\Delta_i\Delta_j G|\leq \lambda(\Delta_i\omega)\mu(\Delta_j\chi),
\end{eqnarray}
where $\lambda$, $\mu$, $\omega$ and $\chi$ are monotone
increasing functions, and
$\Delta_i\omega=\omega(x_i)-\omega(x_{i-1})$,
$\Delta_j\chi=\chi(y_j)-\chi(y_{j-1})$. There are many examples
that the condition (\ref{cond1}) can be checked, e.g.
$f(x,y)=xysin({1\over x}+{1\over y})$, for $\Phi_1(u)=u^p$,
$\Psi_1(u)=u$, where $p>1$. But it seems difficult, if not
impossible, to check Young's condition (\ref{cond2}) for this
example. Needless to say, in the one-parameter case, it is easy to
see that the well-known example of unbounded variation function
$f(x)=xsin{1\over x}$ is of bounded $p$-variation, for any $p>1$.
We can prove multi-parameter Brownian sheet
introduced by Walsh \cite{walsh} in studying stochastic partial
differential equations satisfies definition of $p$,$1$-variation
path ($p>2$), therefore we can define integral w.r.t. Brownian
sheet pathwise and apply this idea to study stochastic PDEs. We
will publish these results in future publications.

We should point out that in this paper we only study the
two-parameter integration of $p, q$-variation path. This is enough
for the purpose of this paper. In this paper, we don't include the
multiplicative integrations as Lyons \cite{terry2},  \cite{terry1}, Lyons and Qian
\cite{terry} investigated for the one-parameter case. We will
study this important problem in future work.

\section{One parameter integral of local time}
\setcounter{equation}{0}

First we recall the definition of $p$-variation
path and its integration theory (see e.g. Young \cite{young1},
Lyons and Qian \cite{terry}).
\begin{defi}\label{def1}
We call a function $f : [x',x'']\to R$ is of bounded $p$-variation
if
\begin{eqnarray}
\sup_E\sum_{i=1}^m|f(x_i)-f(x_{i-1})|^p<\infty,
\end{eqnarray}
where $E:=\{x'=x_0<x_1<\cdots<x_m=x''\}$ is an arbitrary partition of $[x',x'']$. Here $p\geq 1$ is a fixed real number.
\end{defi}

 From Young \cite{young1}, the integral $\int_{x'}^{x''}
f(x)dg(x)=\lim\limits_{m(E)\to 0}\sum\limits_{i=1}^m
f(\xi_i)(g(x_i)-g(x_{i-1}))$ is well defined if $f$ is of bounded
$p$-variation, $g$ is of bounded $q$-variation, and $f$ and $g$
have no common discontinuities. Here $\xi_i\in [x_{i-1},x_i]$,
$p,q\geq 1$, ${1\over p}+{1 \over
q}>1$, $m(E)=\sup\limits_{1\leq i\leq m} (x_i-x_{i-1})$.\\

Consider a continuous semimartingale $X_t$ on a probability space $(\Omega,{\cal F},P)$ with the
decomposition
\begin{eqnarray}
X_t=M_t+V_t,
\end{eqnarray}
where $M_t$ is a local martingale, $V_t$ is an adapted process of bounded variation. Then
there exists semimartingale local time $L_t^x$ of $X_t$ as a nonnegative random field $L=\{L_t^x: (t,x)\in [0,\infty)\times R,\omega\in \Omega\}$.
 Note there is
a different definition of variation established in Revuz and Yor \cite{yor} (see also Marcus and Rosen \cite{rosen}) and the following result is known (P221, Theorem 1.21, \cite{yor}): Let
($\Delta _n$) be a sequence of subdivisions of $[a,b]$ such that $|\Delta _n|\to 0$
as $n\to\infty$, for any nonnegative and finite random variable $S$,
\begin{eqnarray}
\lim\limits _{n\to \infty}\sum _{\Delta _n}(L_S^{a_{i+1}}-L_S^{a_i})^2=4\int _a^bL_S^xdx+\sum \limits
_{a<x\leq b}(L_S^x-L_S^{x-})^2<\infty,
\end{eqnarray}
in probability.
However this variation is not enough to enable us to apply Young's construction of integrals. We need the following new
result to  establish integrations
of local times.

\bigskip

\begin{lem}\label{newlem}
Semimartingale local time $L_t^x$ is of bounded $p$-variation in $x$ for any $t\geq 0$, for any
$p>2$, almost surely.
\end{lem}
{\bf Proof}: By the usual localization argument, we may first
assume that there is a constant $K$ for which $\int_0^t |dV_s|$,
$<M,M>_t\leq K$. By Tanaka's formula
\begin{eqnarray}\label{tanakaf}
L_t^x=(X_t-x)^+-(X_0-x)^+-{\widehat M}_t^x-{\widehat V}_t^x,
\end{eqnarray}
where,
\begin{eqnarray*}
{\widehat M}_t^x=\int_0^t 1_{\{X_s>x\}}dM_s,\  {\widehat V}_t^x=\int_0^t
1_{\{X_s>x\}}dV_s.
\end{eqnarray*}
 First note the function $\varphi_t(x):=(X_t-x)^+-(X_0-x)^+$ is
Lipschitz continuous in $x$ with Lipschitz constant 2, which implies for any
$p>2$ and $a_i<a_{i+1}$
\begin{eqnarray}\label{hz10}
|\varphi_t(a_{i+1})-\varphi_t(a_i)|^p\leq 2^p (a_{i+1}-a_i)^p.
\end{eqnarray}
Secondly, by H\"older inequality, as $V$ is of bounded variation, so
\begin{eqnarray}\label{hz11}
&&|\widehat V_t^{a_{i+1}}-\widehat V_t^{a_i}|^p\nonumber\\
&\leq&|\int_0^t 1_{\{a_i<X_s\leq a_{i+1}\}}|dV_s||^p\nonumber\\
&\leq&c\int_0^t 1_{\{a_i<X_s\leq a_{i+1}\}}|dV_s|,
\end{eqnarray}
where $c$ is a generic constant. To treat ${\widehat M}_t^a$, we
use the method in the proof of Lemma 3.7.5 in Karatzas and Shreve
\cite{ks} or Theorem 6.1.7 in Revuz and Yor \cite{yor},
\begin{eqnarray*}
&&E|\widehat M_t^{a_{i+1}}-\widehat M_t^{a_i}|^p\\
&=&E|\int_0^t 1_{\{a_i<X_s\leq a_{i+1}\}}dM_s|^p\\
&\leq& cE\left(\int_0^t 1_{\{a_i<X_s\leq
a_{i+1}\}}d<M,M>_s\right)^{p\over
2}\\
&=&cE(\int _{a_i}^{a_{i+1}}L_t^xdx)^{p\over 2}\\
&=& c(a_{i+1}-a_{i})^{p\over 2}E({1\over a_{i+1}-a_{i}}
\int _{a_i}^{a_{i+1}}L_t^xdx)^{p\over 2}
\\
&\leq & c(a_{i+1}-a_{i})^{p\over 2}E{1\over a_{i+1}-a_{i}}
\int _{a_i}^{a_{i+1}}(L_t^x)^{p\over 2}dx
\\
&\leq & c(a_{i+1}-a_{i})^{p\over 2}\sup_{x}
E(L_t^x)^{p\over 2}.
\end{eqnarray*}
Here we used Burkholder-Davis-Gundy inequality,
the occupation times formula, Jensen inequality and Fubini theorem. Now from
(\ref{tanakaf}) and using Burkholder-Davis-Gundy inequality again, we have
\begin{eqnarray*}
E(L_t^x)^{p\over 2}&\leq& cE[(X_t-X_0)^{p\over 2}+(\int_0^t
|dV_s|)^{p\over 2}+<M,M>_t^{p\over 4}]
\\
&\leq& cE<M,M>_t^{p\over 4}+cE(\int_0^t |dV_s|)^{p\over
2}+cE<M,M>_t^{p\over 4}<c_1(K,p).
\end{eqnarray*}
Therefore it follows that
\begin{eqnarray}
E|\widehat M_t^{a_{i+1}}-\widehat M_t^{a_i}|^p\leq&
c(a_{i+1}-a_i)^{p\over 2}.
\end{eqnarray}
Here $c$ is a constant depending on $K,p$. Now we use Proposition
4.1.1 in \cite{terry} ($i=1, \gamma>p-1$), for any partition
$\{a_l\}$ of $[a,b]$
\begin{eqnarray*}
\sup_D\sum_l |{\widehat M}_t^{a_{l+1}}-{\widehat M}_t^{a_l}|^p\leq
c(p,\gamma)\sum_{n=1}^{\infty}n^{\gamma}\sum_{k=1}^{2^n}|{\widehat
M}_t^{a_k^n}-{\widehat M}_t^{a_{k-1}^n}|^p.
\end{eqnarray*}
The crucial thing is that the right hand side does not depend on
partition D, where
\begin{eqnarray*}
a_k^n=a+{k\over {2^n}}(b-a),\ k=0,1,\cdots,2^n.
\end{eqnarray*}
We take expectation
\begin{eqnarray*}
&&E\sum_{n=1}^{\infty}n^{\gamma}\sum_{k=1}^{2^n}|{\widehat
M}_t^{a_k^n}-{\widehat M}_t^{a_{k-1}^n}|^p\\
&=&\sum_{n=1}^{\infty}n^{\gamma}\sum_{k=1}^{2^n}E|{\widehat
M}_t^{a_k^n}-{\widehat M}_t^{a_{k-1}^n}|^p\\
&\leq& c\sum_{n=1}^{\infty}n^{\gamma}({{b-a}\over {2^n}})^{{p\over
2}-1}<\infty \ as\  p>2.
\end{eqnarray*}
Therefore
\begin{eqnarray*}
\sum_{n=1}^{\infty}n^{\gamma}\sum_{k=1}^{2^n}|{\widehat
M}_t^{a_k^n}-{\widehat M}_t^{a_{k-1}^n}|^p<\infty\ a.s..
\end{eqnarray*}
It turns out that for any interval $[a,b]\subset R$
\begin{eqnarray*}
\sup_D\sum_l |{\widehat M}_t^{a_{l+1}}-{\widehat M}_t^{a_l}|^p<\infty\
a.s..
\end{eqnarray*}
But we know for each $\omega$, $L_t(a)$ has a compact support in
$a$, say $[-N,N]$ contains its support. Denoting its partition still by
$D:=D_{-N,N}=\{-N=x_0<x_1<\cdots<x_r=N\}$, we obtain
\begin{eqnarray}\label{hz12}
\sup_{D}\sum_i |{\widehat M}_t^{a_{i+1}}-{\widehat M}_t^{a_i}|^p<\infty\
a.s..
\end{eqnarray}
On the other hand, it is easy to see from (\ref{hz10}) that
\begin{eqnarray}\label{hz13}
\sum_i|\varphi_t(a_{i+1})-\varphi_t(a_i)|^p&\leq&
2^p\sum_i(a_{i+1}-a_i)^p\nonumber
\\
&\leq & 2^p[\sum_i(a_{i+1}-a_i)]^p= 2^p(b-a)^p,
\end{eqnarray}
and from (\ref{hz11}) and bounded variation of $V$ that
\begin{eqnarray}\label{hz14}
\sum_i|{\widehat V}_t^{a_{i+1}}-{\widehat V}_t^{a_i}|^p\leq c\int_0^t
1_{\{a<X_s\leq b\}} |dV_s| \leq c\int_0^t
 |dV_s|<\infty.
\end{eqnarray}
Then from (\ref{tanakaf}), (\ref{hz12}), (\ref{hz13}), (\ref{hz14}), we know that
\begin{eqnarray*}
\sup_{D}\sum_i |L_t^{a_{i+1}}-L_t^{a_i}|^p<\infty\ \ \ \ \ \ a.s..
\end{eqnarray*}
Finally we can use the usual localization procedure to  remove the
assumption that $\int_0^t |dV_s|$, $<M,M>_t\leq K$. For this,
define a stopping time for an integer $K>0$: $\tau_K=\inf\{s:
\min\{\int_0^s |dV_r|, <M,M>_s\}> K\}$ if there exists $s$ such
that $\min\{\int_0^s |dV_r|, <M,M>_s\}> K$ and $\tau_K=+\infty$
otherwise. Then the above result shows that there exists $\Omega
_1\subset \Omega$ with $P(\Omega _1)=1$ such that for each
$\omega\in \Omega _1$ and each given integer $K>0$,
\begin{eqnarray*}
\sup_{D}\sum_i |L_{t\wedge \tau_K}^{a_{i+1}}-L_{t\wedge \tau_K}^{a_i}|^p<\infty.
\end{eqnarray*}
Since $\int_0^t |dV_s|(\omega)$ and $<M,M>_t$ are finite almost surely so there exists $\Omega
_2\subset \Omega $ with $P(\Omega _2)=1$ such that for each
$\omega \in \Omega _2$, there exists an integer $K(\omega)>0$ such
that $\int_0^t |dV_s|(\omega)$, $<M,M>_t(\omega)\leq K$. This
leads to  $\tau _K  (\omega)>t$. So for each $\omega \in \Omega
_1\cap \Omega _2$,
\begin{eqnarray*}
\sup_{D}\sum_i |L_t^{a_{i+1}}-L_t^{a_i}|^p<\infty.
\end{eqnarray*}
The result follows as $P(\Omega _1\cap \Omega _2)=1$.
$\hfill\diamond$
\bigskip

Recall the well-known result (see Revuz and Yor \cite{yor}, P220) that for each $t$, the random function $x\to L_t^x$ is a cadlag function hence only admits at most countably many discontinuous points.
Denote $\widehat L_t^x=L_t^x-L_t^{x-}$. Then
\begin{eqnarray}\label{hz15}
\widehat L_t^x=\int _0^t 1_{\{x\}}(X_s)dV_s,
\end{eqnarray}
and for any $a<b$,
\begin{eqnarray}\label{hz31}
\sum \limits _{a< x\leq b}|\widehat L_t^x|=\int _0^t|dV_s|<\infty.
\end{eqnarray}
Write
 \begin{eqnarray}\label{hz102}
 L_t^x=\tilde{L}_t^x+\sum_{x_k^*\leq x} \widehat {L}_t^{x_k^*}.
 \end{eqnarray}
 Here $\tilde{L}_t^x$ is continuous in $x$, and $\{x_k^*\}$ are the
 discontinuous points of $L_t^x$. Denote
\begin{eqnarray}\label{hz103}
h(t,x):=\sum_{x_k^*\leq x} \widehat {L}_t^{x_k^*}.
\end{eqnarray}

\begin{lem}\label{newlem1}
Above defined $h(t,x)$ is of bounded variation in $x$ for each $t$ and of bounded variation in $(t,x)$
for almost every $\omega\in \Omega$.
\end{lem}

\noindent
{\bf Proof}: Let $[-N,N]$ be the support of $L_t(x)$.
 To see that $h(t,x)$ is of locally bounded variation in $x$, consider any partition $D=\{-N=x_0<x_1<\cdots <x_{m-1}<x_m=N\}$, then from (\ref{hz31})
\begin{eqnarray*}
\sum_i|h(t,x_{i+1})-h(t,x_i)|
&=&\sum_i|\sum_{x_i< x_k^*\leq x_{i+1}} \widehat L_t^{x_k^*}|\\
&\leq&\sum_i\sum_{x_i< x_k^*\leq x_{i+1}} |\widehat L_t^{x_k^*}|\\
&=& \sum_{-N< x\leq N}|\widehat L_t^x|<\infty.
\end{eqnarray*}
To see it is of bounded variation in $(t,x)$, consider any partition $D'\times D$, where
 $D'=\{0=t_0<t_1<\cdots <t_{n-1}<t_n=T\}$, $D=\{-N=x_0<x_1<\cdots <x_{m-1}<x_m=N\}$,
 \begin{eqnarray}\label{hz100}
&&\sum_i|h(t_{j+1},x_{i+1})-h(t_{j+1},x_{i})-h(t_{j},x_{i+1})+h(t_j,x_i)|\nonumber\\
&=&\sum_i|\sum_{x_i< x_k^*\leq x_{i+1}} (\widehat
L_{t_{j+1}}^{x_k^*}-\widehat L_{t_{j}}^{x_k^*})|
\nonumber\\
&\leq&\sum_i\sum_{x_i< x_k^*\leq x_{i+1}} |\widehat L_{t_{j+1}}^{x_k^*}-\widehat L_{t_{j}}^{x_k^*}|\nonumber\\
&=& \sum_{-N< x\leq N}|\widehat L_{t_{j+1}}^{x}-\widehat
L_{t_{j}}^{x}|.
\end{eqnarray}
Now applying  (\ref{hz15}) leads to,
\begin{eqnarray}\label{hz101}
&& \sum_{-N< x\leq N}|\widehat L_{t_{j+1}}^{x}-\widehat L_{t_{j}}^{x}|=\sum_{-N< x\leq N}|\int _{t_j}^{t_{j+1}}
 1_{\{x\}}(X_s)dV_s|\nonumber\\
 &
 \leq& \int _{t_j}^{t_{j+1}} 1_{[-N,N]}(X_s)|dV_s|.
\end{eqnarray}
From (\ref{hz100}), (\ref{hz101}) and the bounded variation assumption of $V$, we have
\begin{eqnarray}
&&\sum_j\sum_i|h(t_{j+1},x_{i+1})-h(t_{j+1},x_{i})-h(t_{j},x_{i+1})+h(t_j,x_i)|
\nonumber
\\
&\leq&
\int _{0}^{t} 1_{[-N,N]}(X_s)|dV_s|<\infty.
 \end{eqnarray}
$\hfill\diamond$
\bigskip

Due to the decomposition (\ref{hz102}) of local time, the following integral is therefore defined by
\begin{eqnarray*}
\int _{-\infty}^\infty f(x)d_xL_t^x=\int _{-\infty}^\infty f(x)d_x\tilde{L}_t^x+\int _{-\infty}^\infty f(x)d_xh(t,x).\\
\end{eqnarray*}
The last integral is a Lebesgue-Stieltjes integral, it doesn't matter
whether or not $f$ is continuous as long as it is measurable.
 If $f$ is of finite $p$-variation $(1\leq p<2)$,
we know the integral $\int _{-\infty}^\infty f(x)d_x\tilde{L}_t^x$
 is well defined by Young's integration theory.
\begin{rmk}\label{rem}
 If $f$ belongs to $C^1$, we have
 \begin{eqnarray}\label{fc6}
 \int _{-\infty}^\infty f(x)d_xL_t^x=-\int_{-\infty}^\infty L_t^x df(x).
 \end{eqnarray}
 This is because $L_t^\cdot$ has a compact support for each $t$, so one can always add some points in the partition to
 make $L_t^{x_1}=0$ and $L_t^{x_r}=0$. So
 \begin{eqnarray*}
 &&\int _{-\infty}^\infty f(x)d_x L_t^x\\
 &=&\lim_{m(D)\to 0}\sum _{k=1}^r f(x_{k-1})(L_t^{x_k}-L_t^{x_{k-1}})\\
 &=&\lim_{m(D)\to 0}[\sum _{k=1}^r f(x_{k-1})L_t^{x_k}-\sum _{k=0}^{r-1} f(x_{k})L_t^{x_k}]\\
 &=&-\lim_{m(D)\to 0}\sum _{k=1}^r (f(x_k)-f(x_{k-1}))L_t^{x_k}\\
 &=&-\int_{-\infty}^\infty L_t^x df(x).
 \end{eqnarray*}
 \end{rmk}

Assume $g(x)$ is a left continuous function, we use the standard
regularizing mollifiers to smoothrize $g$ (e.g. see \cite{ks}).
Define
\begin{eqnarray}\label{mollifier}
\rho(x)=\cases {c{\rm e}^{{1\over (x-1)^2-1}}, {\rm \  if } \ x\in (0,2),\cr
0, \ \ \ \ \ \ \ \ \ \ \ \  {\rm otherwise.}}
\end{eqnarray}
Here $c$ is chosen such that $\int _0^2\rho(x)dx=1$.
Take $\rho_n(x)=n\rho(nx)$ as mollifiers. Define
\begin{eqnarray*}
g_n(x)=\int _{-\infty}^{+\infty}\rho_n(x-y)g(y)dy, \ \ n\geq 1.
\end{eqnarray*}
Then $g_n(x)$ is smooth and
\begin{eqnarray}\label{f1}
g_n(x)=\int _0^2\rho(z)g(x-{z\over n})dz, \ \ n\geq 1.
\end{eqnarray}
Using Lebesgue's dominated convergence theorem, one can prove that
as $n\to \infty$, $g_n(x)\to  g(x)$.

\begin{thm}\label{fcr1}
Let $g(x)$ be a left continuous function with finite $p$-variation
in $x$, $1\leq p<2$, $g_n(x)$ be defined in (\ref{f1}). Then
\begin{eqnarray}\label{con2}
\int_{-\infty}^{\infty}g_n(x)d_x\tilde{L}_t^x\to
\int_{-\infty}^{\infty}g(x)d_x\tilde{L}_t^x , \  as \ n\to \infty.
\end{eqnarray}
\end{thm}
{\bf Proof}: Let $\delta >0$ satisfy ${1\over {2+\delta}}+{1\over
p}>1$. From Lemma \ref{newlem}, $\tilde{L}_t^x$ is of bounded
($2+\delta$)-variation in $x$ uniformly in $t$. From \cite{young},
$g(x)$ being of bounded $p$-variation, $1\leq p<2$, is equivalent
to that for any partition $D:=D_{-N,N}=
\{-N={x_0}<{x_1}<\cdots<{x_r}=N\}$ defined as before, there is an
increasing function $w$ such that
\begin{eqnarray*}
&&|g(x_{l+1})-g(x_l)|\leq (w(x_{l+1})-w(x_l))^{1\over p}, \
\forall x_{l}, x_{l+1}\in D,
\end{eqnarray*}
where $w(x)$ is the total $p$-variation of $f$ in the interval
$[-N-2, x]$. Using H\"older inequality, we get
\begin{eqnarray*}
&&\sup_D \sum_{l=1}^r|g_n(x_l)-g_n(x_{l-1})|^p\\
&=&\sup_D \sum_{l=1}^r \left|\int_0^2\rho(z)[g(x_l-{z\over n})-g(x_{l-1}-{z\over n})]dz\right|^p\\
&\leq&M_1\sup_D \sum_{l=1}^r\left(\int_0^2|g(x_l-{z\over n})-g(x_{l-1}-{z\over n})|^p dz\right)\\
&\leq&M_1\int_0^2\sup_D \sum_{l=1}^r |g(x_l-{z\over n})-g(x_{l-1}-{z\over n})|^p dz\\
&\leq&M_{1}\int_0^2  (w(N-{z\over n})-w(-N-{z\over n}))dz,
\end{eqnarray*}
where $M_1$ is a constant. As
\begin{eqnarray*}
w(N-{z\over n})-w(-N-{z\over n})\leq w(N),
\end{eqnarray*}
so
\begin{eqnarray}
\sup_D \sum_{l=1}^r|g_n(x_l)-g_n(x_{l-1})|^p\leq 2M_1w(N)<\infty,
\end{eqnarray}
which means that $g_n(x)$ is of bounded $p$-variation in $x$
uniformly in $n$. Then (\ref{con2}) follows from Young's
(\cite{young1} or \cite{young}) convergence theorem we can get the
result directly. $\hfill\diamond$ \vskip5pt
\begin{rmk}
From the Lebesgue's dominated convergence theorem, for $g$ in the
above theorem, we know
\begin{eqnarray*}
\int_{-\infty}^\infty g_n(x)d_xh(t,x)\to \int_{-\infty}^\infty
g(x)d_xh(t,x),\ as\ n \to \infty.
\end{eqnarray*}
With Theorem \ref{fcr1}, it follows that
\begin{eqnarray}\label{add1}
\int_{-\infty}^\infty g_n(x)d_xL_t^x\to \int_{-\infty}^\infty
g(x)d_x{L}_t^x,\ as\ n \to \infty.
\end{eqnarray}
\end{rmk}

 Using the above theorem, we can get an extension of It$\hat {\rm o}$'s Formula.

  \begin{thm}\label{fcr11}
 Let $X=(X_t)_{t\geq0}$ be a continuous semi-martingale and $f: R\to R$ be a left continuous, locally bounded function and have left derivative $\nabla ^-f(x)$ being left continuous and locally bounded. Assume $\nabla ^-f(x)$ is of bounded $q$-variation, where $1\leq q <2$. Then we have
the following change-of-variable formula
 \begin{eqnarray}
 f(X_t)=f(X_0)+\int_0^t \nabla ^-f(X_s)dX_s-\int_{-\infty}^{\infty}\nabla ^-f(x)d_x L_t^x,
 \end{eqnarray}
 where $L_t^x$ is the local time of $X_t$ at $x$.
 \end{thm}
 {\bf Proof}: The integral $\int_{-\infty}^\infty\nabla ^-f(x)d_x L_t^x$ is defined pathwise as a combination of rough path integral and
 Lebesgue-Stieltjes integral.
 We may quote the proof in \cite{ks} and define
\begin{eqnarray*}
f_n(x)=\int _{-\infty}^{+\infty}\rho_n(x-y)f(y)dy, \ \ n\geq 1.
\end{eqnarray*}
The convergence of all terms except the second order derivative
term are the same as in the proof in \cite{ks}. By occupation
times formula and Remark \ref{rem}, the second order derivative
term is
 \begin{eqnarray*}
 {1\over2}\int_0^t \Delta f_n(X_s)d<M>_s
&=&\int_{-\infty}^\infty \Delta f_n(x) L_t^x dx\\
&=&\int_{-\infty}^\infty L_t^x d\nabla f_n(x)\\
&=&-\int_{-\infty}^\infty \nabla f_n(x) d_xL_t^x.
\end{eqnarray*}
It follows from (\ref{add1}) that,
 \begin{eqnarray*}
 &&{1\over2}\int_0^t \Delta f_n(X_s)d<M>_s\to -\int_{-\infty}^\infty \nabla ^-f(x) d_xL_t^x,
 \end{eqnarray*}
 when $n\to \infty$. Our claim is asserted.
$\hfill\diamond$
\vskip5pt

Needless to say, there are many cases that Theorem \ref{fcr11} works, but other extensions of
It${\hat {\rm o}}$'s formula do not apply immediately. The following is an obvious example:

\begin{exap} Consider a function $f(x)=x^3\cos {1\over x}$ for $x\neq 0$ and $f(0)=0$.
This function is $C^1$ and its derivative is
$f^{\prime}(x)=3x^2\cos {1\over x}+x\sin {1\over x}$ for $x\neq 0$
and $f^{\prime}(0)=0$. It is easy to see that $f^{\prime}$ is not
of bounded variation, but of $p$-variation for any $p>1$ (see
Example 3.1 for a proof in a more complicated case). So Theorem
\ref{fcr11} can be used, while Meyer's formula cannot apply to
this situation.
\end{exap}
\section{Two-parameter $p, q$-variation path integrals }
\setcounter{equation}{0}

      In this section, the following notations are used:
$\Phi(u)$, $\Psi(u)$, $\Phi_1(u)$, $\Psi_1(u)$ denote continuous
functions strictly increasing from $0$ to $\infty$ with u, where
$u\geq 0$ is a variable, and $\Phi(0)$, $\Psi(0)$, $\Phi_1(0)$,
$\Psi_1(0)\equiv0$; $\varphi(u)$, $\psi(u)$, $\varphi_1(u)$,
$\psi_1(u)$ denote the inverse functions of  $\Phi(u)$, $\Psi(u)$,
$\Phi_1(u)$, $\Psi_1(u)$, respectively; $\omega$, $\chi$ are
monotone increasing functions of one variable, $x$ or $y$.

      Firstly in the following we will define a two-parameter Young integral
      $\int_{x'}^{x''}\int_{y'}^{y''}F(x,y)d_{x,y}G(x,y)$. We will use some idea
      from Young \cite{young}. But Young's condition is very strong and the class of
      functions that satisfy Young's condition is restricted. In particular,
      Young's condition does not seem to include the class of functions of bounded variation
      and many important examples. We give a new and weaker condition for the integration in
      this section. We will use Lyons' idea of control functions to simplify our proof.
      One can see our condition is a natural extension of locally bounded multi-dimensional
      L-S measure.
First, if $F(x,y)$ is a simple function, say
\begin{eqnarray*}
F(x,y)=\sum_{i=1}^M\sum_{j=1}^{M'}F(x_{i-1},y_{j-1})1_{\{x_{i-1}< x \leq x_i, y_{j-1}< y \leq y_j\}},
\end{eqnarray*}
as normal we can see that the integral of the simple function can be defined as
\begin{eqnarray*}
&&\int_{x'}^{x''}\int_{y'}^{y''}F(x,y)d_{x,y}G(x,y)\\
&=&\sum_{i=1}^M\sum_{j=1}^{M'}F(x_{i-1},y_{j-1})\left(G(x_i,y_j)-G(x_{i-1},y_j)-G(x_i,y_{j-1})+G(x_{i-1},y_{j-1})\right).
\end{eqnarray*}

Before we proceed, we need the following definition.

\begin{defi} Let $E\times E':=\{x'=x_0<x_1<\cdots<x_N=x'', y'=y_0<y_1<\cdots<y_{N'}=y''\}$ be an arbitrary partition of $[x',x'']\times[y',y'']$.
We call $F(x,y)$ is of bounded $\Phi$-variation in $x$ uniformly in $y$, if
\begin{eqnarray*}
\sup\limits_{y\in[y',y'']}\sup\limits_E\sum\limits_{k=1}^N \Phi(|F(x_k,y)-F(x_{k-1},y)|)< \infty.
\end{eqnarray*}
We call $G(x,y)$ is of bounded $\Phi_1$, $\Psi_1$-variation in ($x, y$), if
\begin{eqnarray}
\sup_{E\times E'}\sum_{j=1}^{N'}\Psi_1\left(\sum_{i=1}^N\Phi_1(|\Delta_i\Delta_j G|)\right)<\infty,
\end{eqnarray}
where
\begin{eqnarray*}
\Delta_i\Delta_j G:=G(x_i,y_j)-G(x_{i-1},y_j)-G(x_i,y_{j-1})+G(x_{i-1},y_{j-1}).
\end{eqnarray*}
If $\Psi_1(u)=u$, we call $G(x,y)$ is of bounded
$\Phi_1$-variation in $(x,y)$.
 If $\Phi_1(u)=u^p$, $\Psi_1(u)=u^q$, $p,q\geq 1$, we call
$G(x,y)$ is of bounded $p,q$-variation in $(x,y)$.
\end{defi}

 In the following, we will give an example of $p,1$-variation ($p>1$)
 function.
 \begin{exap}
Consider
\begin{eqnarray*}
 f(x,y)=xysin({1\over x}+{1\over y}), \ 0<x,y\leq 1,\
f(0,y)=f(x,0)=f(0,0)=0.
\end{eqnarray*}
This is a continuous function of unbounded variation but of
bounded $p,1$-variation ($p>1$).
 To see it is of unbounded variation, we take the partition
  $E_1\times E_2=\{0<{1\over {n\pi +{\pi\over 2} -1}}<{1\over {n\pi -1}}<\cdots<{1\over {\pi-1}}<1, 0<1\}$,
 \begin{eqnarray*}
 &&\sum_{i,j}|\Delta_i\Delta_jf|=\sum_i|x_isin({1\over{x_i}}+1)-x_{i-1}sin({1\over {x_{i-1}}}+1)|\\
 &=&\sum_i|{1\over {i \pi +{\pi\over 2} -1}}sin(i\pi+{\pi\over 2})-{1\over {i\pi -1}}sin(i\pi)|\\
 &=&\sum_i{1\over {i \pi +{\pi\over 2} -1}}\\
 &=&\infty.
 \end{eqnarray*}
To see it is of bounded $p,1$-variation for any $p>1$, consider
any partition $E\times E'$
\begin{eqnarray}\label{hz32}
&&\sum_{i,j}|\Delta_i\Delta_jf|^p\nonumber\\
&=&\sum_{i,j}\bigg|x_iy_jsin({1\over {x_i}}+{1\over {y_j}})-x_{i-1}y_jsin({1\over {x_{i-1}}}+{1\over {y_j}})\nonumber\\
&&\hskip 0.7cm-x_iy_{j-1}sin({1\over {x_i}}+{1\over {y_{j-1}}})+x_{i-1}y_{j-1}sin({1\over {x_{i-1}}}+{1\over {y_{j-1}}})\bigg|^p\nonumber\\
&=&\sum_{i,j}\bigg|y_j\Big[x_isin({1\over {x_i}}+{1\over {y_j}})-x_{i-1}sin({1\over {x_{i-1}}}+{1\over {y_j}})\Big]\nonumber\\
&&\hskip 0.7cm -y_{j-1}\Big[x_isin({1\over {x_i}}+{1\over {y_{j-1}}})-x_{i-1}sin({1\over {x_{i-1}}}+{1\over {y_{j-1}}})\Big]\bigg|^p\nonumber\\
&=&\sum_{i,j}\bigg| y_j(x_i-x_{i-1})sin({1\over {x_i}}+{1\over {y_j}})+y_jx_{i-1}\Big[sin({1\over {x_i}}+{1\over {y_j}})-sin({1\over {x_{i-1}}}+{1\over {y_j}})\Big]\nonumber\\
&&\hskip 0.7cm -y_{j-1}(x_i-x_{i-1})sin({1\over {x_i}}+{1\over {y_{j-1}}})\nonumber\\
&&\hskip 0.7cm -y_{j-1}x_{i-1}\Big[sin({1\over {x_i}}+{1\over {y_{j-1}}})-sin({1\over {x_{i-1}}}+{1\over {y_{j-1}}})\Big]\bigg|^p\nonumber\\
&=&\sum_{i,j}\bigg|(x_i-x_{i-1})\Big[y_jsin({1\over {x_i}}+{1\over {y_j}})-y_{j-1}sin({1\over {x_i}}+{1\over {y_{j-1}}})\Big]\nonumber\\
&&\hskip 0.7cm +x_{i-1}\Big[y_j\big((sin({1\over {x_i}}+{1\over {y_j}})-sin({1\over {x_{i-1}}}+{1\over {y_j}})\big)\nonumber\\
&&\hskip 1.9cm -y_{j-1}\big(sin({1\over {x_i}}+{1\over {y_{j-1}}})-sin({1\over {x_{i-1}}}+{1\over {y_{j-1}}})\big)\Big]\bigg|^p\nonumber\\
&=&\sum_{i,j}\bigg|(x_i-x_{i-1})\Big[(y_j-y_{j-1})sin({1\over {x_i}}+{1\over {y_j}})\nonumber\\
&& \hskip1.9cm
+y_{j-1}\big(sin({1\over {x_i}}+{1\over {y_{j}}})-sin({1\over {x_i}}+{1\over {y_{j-1}}})\big)\Big]\nonumber\\
&&\hskip 0.7cm +x_{i-1}\Big[(y_j-y_{j-1})\big(sin({1\over {x_i}}+{1\over {y_j}})-sin({1\over {x_{i-1}}}+{1\over {y_j}})\big)\nonumber\\
&&\hskip 1.9cm +y_{j-1}\big(sin({1\over {x_i}}+{1\over {y_j}})-sin({1\over {x_{i-1}}}+{1\over {y_j}})\nonumber\\
&&\hskip 3cm-sin({1\over {x_i}}+{1\over {y_{j-1}}})+sin({1\over {x_{i-1}}}+{1\over {y_{j-1}}})\big)\Big]\bigg|^p\nonumber\\
&\leq& c_p\bigg\{\sum_{i,j}\Big|(x_i-x_{i-1})(y_j-y_{j-1})sin({1\over {x_i}}+{1\over {y_j}})\Big|^p\nonumber\\
&&\hskip 0.7cm+\sum_{i,j}\Big|y_{j-1}(x_i-x_{i-1})\big(sin({1\over {x_i}}+{1\over {y_{j}}})-sin({1\over {x_i}}+{1\over {y_{j-1}}})\big)\Big|^p\nonumber\\
&&\hskip 0.7cm+\sum_{i,j}\Big|x_{i-1}(y_j-y_{j-1})\big(sin({1\over {x_i}}+{1\over {y_j}})-sin({1\over {x_{i-1}}}+{1\over {y_j}})\big)\Big|^p\nonumber\\
&&\hskip 0.7cm+\sum_{i,j}\Big|x_{i-1}y_{j-1}\big(sin({1\over {x_i}}+{1\over {y_j}})-sin({1\over {x_{i-1}}}+{1\over {y_j}})\nonumber\\
&&\hskip 3.5cm-sin({1\over {x_i}}+{1\over {y_{j-1}}})+sin({1\over {x_{i-1}}}+{1\over {y_{j-1}}})\big)\Big|^p\bigg\}\nonumber\\
&&:=c_p(I+II+III+IV),
\end{eqnarray}
where $c_p$ is a constant. It's easy to see that
\begin{eqnarray}\label{hz33}
I\leq \sum_{i,j}(x_i-x_{i-1})^p(y_j-y_{j-1})^p\leq 1.
\end{eqnarray}
For $II$, as $|sinx|\leq x$, so
\begin{eqnarray*}
&&II\leq 2^{p-1}\sum_{i,j}{y^p_{j-1}}(x_i-x_{i-1})^p\Big|sin({1\over {x_i}}+{1\over {y_{j}}})-sin({1\over {x_i}}+{1\over {y_{j-1}}})\Big|\\
&&\hskip 0.5cm=2^{p-1}\sum_{i,j}{y^p_{j-1}}(x_i-x_{i-1})^p\cdot\Big|2cos{{2\over {x_i}}+{1\over {y_{j}}}+{1\over {y_{j-1}}}\over 2}sin{{{1\over {y_{j}}}-{1\over {y_{j-1}}}}\over 2}\Big|\\
&&\hskip 0.5cm\leq
2^{p-1}\sum_{i,j}{y^p_{j-1}}(x_i-x_{i-1})^p\left({1\over{y_{j-1}}}-{1\over{y_{j}}}\right)\\
&&\hskip 0.5cm=
2^{p-1}\sum_{i}(x_i-x_{i-1})^p\sum_j{y^p_{j-1}}\left({1\over{y_{j-1}}}-{1\over{y_{j}}}\right).
\end{eqnarray*}
It is obvious that
\begin{eqnarray*}
\sum_{i}(x_i-x_{i-1})^p<\infty.
\end{eqnarray*}
And also because
\begin{eqnarray}\label{bb1}
\sum_j{y^p_{j-1}}\left({1\over{y_{j-1}}}-{1\over{y_{j}}}\right)
&\leq&\sum_j{y^p_{j-1}}{{y_j-y_{j-1}}\over{y^2_{j-1}}}\nonumber\\
&=&\sum_j{y^{p-2}_{j-1}}(y_j-y_{j-1})\nonumber\\
&\leq &\int_0^1 y^{p-2} dy={1\over {p-1}}.
\end{eqnarray}
So we get $II<\infty$.  Similar to the discussion of $II$, we can
also prove that $III<\infty$. About $IV$,
\begin{eqnarray*}
&&IV\\
&=&\sum_{i,j}{x_{i-1}^p}{y_{j-1}^p}\Big|2cos{{{1\over {x_i}}+{1\over {x_{i-1}}}+{2\over {y_j}}}\over 2}sin{{{1\over {x_i}}-{1\over {x_{i-1}}}}\over 2}\\
&&\hskip 2.5cm-2cos{{{1\over {x_i}}+{1\over {x_{i-1}}}+{2\over {y_{j-1}}}}\over 2}sin{{{1\over {x_i}}-{1\over {x_{i-1}}}}\over 2}\Big|^p\\
&=&2^p\sum_{i,j}{x_{i-1}^p}{y_{j-1}^p}\Big|sin{{{1\over {x_i}}-{1\over {x_{i-1}}}}\over 2}\Big|^p\cdot\Big|cos{{{1\over {x_i}}+{1\over {x_{i-1}}}+{2\over {y_j}}}\over 2}-cos{{{1\over {x_i}}+{1\over {x_{i-1}}}+{2\over {y_{j-1}}}}\over 2}\Big|^p\\
&=&2^p\sum_{i,j}{x_{i-1}^p}{y_{j-1}^p}\Big|sin{{{1\over
{x_i}}-{1\over {x_{i-1}}}}\over 2}\Big|^p\\
&&
\hskip2.5cm \cdot\Big|-2sin{{{1\over
{x_i}}+{1\over {x_{i-1}}}+{1\over {y_j}}+{1\over {y_{j-1}}}}\over
2}sin{{{1\over {y_j}}-{1\over {y_{j-1}}}}\over 2}\Big|^p\\
&\leq& 2^p\cdot
2^{p}\sum_{i,j}{x_{i-1}^p}{y_{j-1}^p}{1\over2}({1\over
{x_{i-1}}}-{1\over {x_i}})\cdot{1\over2}({1\over
{y_{j-1}}}-{1\over {y_j}})\\
&=&2^{2p-2}\sum_i{x_{i-1}^p}({1\over {x_{i-1}}}-{1\over
{x_i}})\sum_j{y_{j-1}^p}({1\over {y_{j-1}}}-{1\over {y_j}})\\
&\leq& {2^{2p-2}\over (p-1)^2} \ ,
\end{eqnarray*}
following from a similar argument as in (\ref{bb1}). So the function $f(x,y)=xysin({1\over x}+{1\over
y})$, $0<x,y\leq 1$, $f(0,y)=f(x,0)=f(0,0)=0$, is of bounded
$p,1$-variation for any $p>1$. Moreover, from the above proof, we
can see for this function $f(x,y)$ on $(x,y)\in [0,\delta_1]\times [0,\delta_2]$, its $p,1$-variation tends to $0$ when either
$\delta_1$ or $\delta_2$ decreases to $0$. $\hfill\diamond$
 \vskip5pt
\end{exap}

We say a function $f(x,y)$ has a jump at $(x_1,y_1)$ if there
exists an $\varepsilon >0$ such that for any $\delta>0$, there exists
$(x_2,y_2)$ satisfying  $\max\{|x_1-x_2|,|y_1-y_2|\}<\delta$ and
$|f(x_2,y_2)-f(x_1,y_2)-f(x_2,y_1)+f(x_1,y_1)|>\varepsilon$. For a
function $G(x,y)$ of bounded $\Phi_1, \Psi_1$-variation, for any
given $\varepsilon>0$, it is easy to see that there exists a $\delta (\varepsilon)>0$ and
 a finite number of jump points $\{(x_1,y_1),\cdots,(x_{n_0},y_{m_0})\}$ such that
  $|G(x,y)-G(x,\tilde y)-G(\tilde x,y)+G(\tilde x,\tilde y)|<\varepsilon$ whenever
$\max\{|\tilde x-x|,|\tilde y-y|\}<\delta(\varepsilon)$,
 $[\tilde x,x]\cap\{x_1,\cdots,x_{n_0}\}= \emptyset$ and
 $ [\tilde y,y]\cap\{y_1,\cdots,y_{m_0}\}= \emptyset$.
Denote $H_0\times
H'_0:=\{x_1,\cdots,x_{n_0}\}\times\{y_1,\cdots,y_{m_0}\}$.
In the following, we assume the following finite large jump
condition: for any $\varepsilon >0$, there exists at
most finite many points $\{x_1, \cdots ,x_{n_1}\}$, $\{y_1, \cdots ,y_{m_1}\}$
and a constant $\delta (\varepsilon)
>0$ such that the total $\Phi_1,
\Psi_1$-variation of $G$ on $[x, x+\delta]\times [y',y'']$ is smaller than $\varepsilon$
if $[x, x+\delta]\cap \{x_1, \cdots ,x_{n_1}\}= \emptyset$, and
the total $\Phi_1,
\Psi_1$-variation of $G$ on $[x', x'']\times [y,y+\delta]$ is smaller than $\varepsilon$
if $[y, y+\delta]\cap \{y_1, \cdots ,y_{m_1}\}= \emptyset$.
Denote $H\times
H':=\{x_1,\cdots,x_{n_1}\}\times\{y_1,\cdots,y_{m_1}\}$. It is
obviously that $H\times H'\supset H_0\times H'_0$.
There are many examples of bounded $\Phi_1,\Psi_1$-variation functions
that satisfy the finite large jump condition. But it is not clear whether or not
the bounded $\Phi_1,\Psi_1$-variation condition implies automatically the finite large
jump condition in the two parameter case although this is true in the
one parameter case.

   Denote by $\omega(x_k)$ the total uniform $\Phi$-variation of
$F$ in $x$ in the interval $[x',x_k]$ and $\chi(y_l)$ the total
uniform $\Psi$-variation of $F$ in $y$ in the interval $[y',y_l]$.
For the partition $E\times E'$, denote by $m(E\times E')$ the mesh
of the partition.

\vskip5pt

We need the following simple inequalities: Let $f$ be a nonnegative and nondecreasing function,
then
\begin{eqnarray}\label{hz38}
\sum_{p=0}^{\infty} 2^{p-1}f({1\over 2^p})\leq \sum _{m=1}^{\infty}f({1\over m})\leq
\sum_{p=0}^{\infty} 2^{p}f({1\over 2^p}),
\end{eqnarray}
and for any $v\geq 1$,
\begin{eqnarray}\label{hz39}
\sum_{p=v}^{\infty} 2^{p-1}f({1\over 2^p})\leq \sum
_{m=2^{v-1}+1}^{\infty}f({1\over m})\leq \sum_{p=v-1}^{\infty}
2^{p}f({1\over 2^{p}}),
\end{eqnarray}
if the series $\sum _{m=1}^{\infty}f({1\over m})$ is convergent.
These inequalities were also used in the proof of Young's main results.
 We listed them here only for the purpose  to make
the proof of the following theorem easier to understand. The proof
is elementary and omitted.

\begin{thm}\label{ffcr1}
Let $F(x,y)$ be a continuous function of bounded $\Phi$-variation
in $x$ uniformly in $y$, and be of bounded $\Psi$-variation in $y$
uniformly in $x$; $G(x,y)$ be of bounded $\Phi_1,
\Psi_1$-variation in ($x,y$) and satisfy the finite large jump condition,
where $\Psi_1$, $\Phi_1$ are convex
functions. If there exist monotone increasing functions $\varrho$
and $\sigma$ subject to $\varrho(u)\sigma(u)=u$ such that
\begin{eqnarray}
\sum\limits_{m,n}\varrho[\varphi({1\over n})]\sigma[\psi({1\over m})]\varphi_1[{1\over n}\psi_1({1\over m})]<\infty,
\end{eqnarray}
 then the integral
\begin{eqnarray}\label{chunrong1}
&&\int_{x'}^{x''}\int_{y'}^{y''}F(x,y)d_{x,y}G(x,y)\nonumber\\
&=&\lim_{m(E\times E')\to
0}\sum_{i=1}^N\sum_{j=1}^{N'}F(x_{i-1},y_{j-1})\Delta_i\Delta_j G
\end{eqnarray}
 is well defined with all partitions $E\times E'$ of
$[x',x'']\times [y',y'']$ including suitable finite sets $H\times H'$ as defined above, i.e. for any given
$\varepsilon>0$, we can determine finite sets $H$ and $H^{\prime}$ of variables $x$ and $y$ respectively such that
\begin{eqnarray*}
|\int_{x'}^{x''}\int_{y'}^{y''}F(x,y)d_{x,y}G(x,y)
-\sum_{i=1}^N\sum_{j=1}^{N'}F(x_{i-1},y_{j-1})\Delta_i\Delta_j G|<\varepsilon
\end{eqnarray*}
as long as the partition $E=\{x^{\prime}=x_1<x_2<\cdots<x_N=x^{\prime\prime}\}$ and
$E^{\prime}=\{y^{\prime}=y_1<y_2<\cdots<y_{N^{\prime}}=y^{\prime\prime}\}$
includes $H$ and $H^{\prime}$ among their points of divisions respectively.

\end{thm}
{\bf Proof}: For any partition $E\times
E':=\{x'=x_0<x_1<\cdots<x_N=x'', y'=y_0<y_1<\cdots<y_{N'}=y''\}$,
consider
\begin{eqnarray*}
F_{E,E'}(x,y):=\sum_{i=1}^N\sum_{j=1}^{N'}F(x_{i-1},y_{j-1})1_{\{x_{i-1}\leq
x < x_i, y_{j-1}\leq y < y_j\}},
\end{eqnarray*}
then
\begin{eqnarray*}
S(E,E')&:=&S_{F_{E,E'}}(E,E'):=\int_{x'}^{x''}\int_{y'}^{y''}F_{E,E'}(x,y)d_{x,y}G(x,y)\\
&=&\sum_{i=1}^N\sum_{j=1}^{N'}F(x_{i-1},y_{j-1})\Delta_i\Delta_j
G.
\end{eqnarray*}
From the assumption of $F$, let $P$, $Q$ be the total $\Phi$- and $\Psi$-variation of $F$ in $x$ and $y$ respectively, so
\begin{eqnarray*}
\sum\limits_{k=1}^N \Phi(|F(x_k,y)-F(x_{k-1},y)|)\leq P,\\
\sum\limits_{l=1}^{N'} \Psi(|F(x,y_l)-F(x,y_{l-1})|)\leq Q,
\end{eqnarray*}
which are equivalent to
\begin{eqnarray*}
|F(x_k,y)-F(x_{k-1},y)|\leq \varphi(\omega(x_k)-\omega(x_{k-1})),\ k=1,2,\cdots,N,\\
|F(x,y_l)-F(x,y_{l-1})|\leq \psi(\chi(y_l)-\chi(y_{l-1})),\ l=1,2,\cdots,N'.
\end{eqnarray*}
Obviously, if $y_{j-1}\leq y< y_{j}$, $j=1,\cdots,N'$,
\begin{eqnarray*}
&&|F_{E,E'} (x_k,y)-F_{E,E'} (x_{k-1},y)|\\
&=&|F(x_k,y_{j-1})-F(x_{k-1},y_{j-1})|\\
&\leq&\varphi(\omega(x_k)-\omega(x_{k-1})),\ \ k=1,2,\cdots,N,
\end{eqnarray*}
and if $x_{i-1}\leq x< x_{i}$, $i=1,\cdots,N$,
\begin{eqnarray*}
&&|F_{E,E'} (x,y_l)-F_{E,E'}(x,y_{l-1})|\\
&=& |F(x_{i-1},y_l)-F(x_{i-1},y_{l-1})|\\
&\leq&\psi(\chi(y_l)-\chi(y_{l-1})),\ \ l=1,2,\cdots,N'.
\end{eqnarray*}
Because $\omega$ and $\chi$ are both increasing functions, so we can define a sequence of finite sets:
\begin{eqnarray*}
&&E_p:\ \omega(x^{(2)})-\omega(x^{(0)})\leq 2^{-p}P,\\
&&E'_q:\ \chi(y^{(2)})-\chi(y^{(0)}) \leq 2^{-q}Q,
\end{eqnarray*}
where $x^{(0)}, x^{(2)}$ are two consecutive points of $E_p$,
$y^{(0)}, y^{(2)}$ are consecutive points of $E'_q$. Such $E_p$
and $E'_q$ can be defined by induction as for each $p\geq 0$ on
$E_p$ either $\omega(x^{(2)})-\omega(x^{(0)})\leq 2^{-(p+1)}P$ or
$2^{-(p+1)}P<\omega(x^{(2)})-\omega(x^{(0)})\leq 2^{-p}P$. In the
latter case, we insert one point between such a pair to get
$E_{p+1}$ (we insert at most $2^p$ points) such that
\begin{eqnarray*}
E_{p+1}:\ \omega(x^{(2)})-\omega(x^{(1)})\leq 2^{-(p+1)}P,
\end{eqnarray*}
where $x^{(1)}$, $x^{(2)}$ are consecutive points of $E_{p+1}$. In
the same way, we can get $E'_{q+1}$ such that
\begin{eqnarray*}
E'_{q+1}:\ \chi(y^{(2)})-\chi(y^{(1)}) \leq 2^{-(q+1)}Q,
\end{eqnarray*}
where $y^{(1)}, y^{(2)}$ are consecutive points of $E'_{q+1}$. In
$E_{p+1}$, there are at most $2^{p+1}$ points and in $E'_{q+1}$,
there are at most $2^{q+1}$ points. We will prove our theorem in
four steps. \vskip 5pt
{\bf Step 1}: Note
\begin{eqnarray}\label{equ1}
&&S(E_{p+1}, E'_{q+1})-S(E_{p}, E'_{q+1}) -S(E_{p+1}, E'_{q}) +S(E_{p}, E'_{q})\nonumber\\
&=&\sum\limits_{i=1,3,5\cdots,2^{p+1}-1} \sum\limits_{j=1,3,5\cdots,2^{q+1}-1}\bigg[F(x_{i-1},y_{j-1})\Delta_i\Delta_j G \nonumber\\
&&
+F(x_{i-1},y_{j})\Delta_i\Delta_{j+1} G+F(x_{i},y_{j})\Delta_{i+1}\Delta_{j+1} G+F(x_{i},y_{j-1})\Delta_{i+1}\Delta_j G\nonumber\\
&&
-F(x_{i-1},y_{j-1})\bigg(G(x_{i+1},y_{j})-G(x_{i-1},y_{j})\nonumber\\
&&\hskip4.5cm -G(x_{i+1},y_{j-1})+G(x_{i-1},y_{j-1})\bigg)\nonumber\\
&& -F(x_{i-1},y_{j})\bigg(G(x_{i+1},y_{j+1})-G(x_{i+1},y_{j})\nonumber\\
&&\hskip4.5cm -G(x_{i-1},y_{j+1})+G(x_{i-1},y_{j})\bigg)\nonumber\\
&& -F(x_{i-1},y_{j-1})\bigg(G(x_{i},y_{j+1})-G(x_{i-1},y_{j+1})\nonumber\\
&&\hskip4.5cm -G(x_{i},y_{j-1})+G(x_{i-1},y_{j-1})\bigg)\nonumber\\
&&-F(x_{i},y_{j-1})\bigg(G(x_{i+1},y_{j+1})-G(x_{i+1},y_{j-1})-G(x_{i},y_{j+1})+G(x_{i},y_{j-1})\bigg)\nonumber\\
&&+F(x_{i-1},y_{j-1})\bigg(G(x_{i+1},y_{j+1})-G(x_{i+1},y_{j-1})
\nonumber\\
&&\hskip4.5cm -G(x_{i-1},y_{j+1})+G(x_{i-1},y_{j-1})\bigg)\bigg]\nonumber\\
&=&\sum\limits_{i=1,3,5\cdots,2^{p+1}-1}
\sum\limits_{j=1,3,5\cdots,2^{q+1}-1}\Big(\Delta_i\Delta_jF\Big)\Big(\Delta_{i+1}\Delta_{j+1}G\Big).
\end{eqnarray}
Because
\begin{eqnarray*}
|\Delta_i\Delta_jF|&\leq& |F(x_{i},y_{j})-F(x_{i-1},y_{j})|+|F(x_{i},y_{j-1})-F(x_{i-1},y_{j-1})|\\
&\leq& 2\varphi (2^{-(p+1)}P)\leq 2\varphi(2^{-p}P),
\end{eqnarray*}
and also
\begin{eqnarray*}
|\Delta_i\Delta_jF|&\leq& |F(x_{i},y_{j})-F(x_{i},y_{j-1})|+|F(x_{i-1},y_{j})-F(x_{i-1},y_{j-1})|\\
&\leq& 2\psi(2^{-(q+1)}Q)\leq 2\psi(2^{-q}Q),
\end{eqnarray*}
so it is easy to see
\begin{eqnarray}\label{equ2}
|\Delta_i\Delta_jF|&\leq& 2\varrho[\varphi (2^{-p}P)]\sigma[\psi(2^{-q}Q)]
\end{eqnarray}
for some increasing functions $\varrho$, $\sigma$ satisfying $\varrho(u)\sigma(u)=u$.

For the function $G$,  let $M$ be its total
$\Phi_1,\Psi_1$-variation, then
\begin{eqnarray*}
\sum_{j=1}^{2^q}\Psi_1\left(\sum_{i=1}^{2^p} \Phi_1 (|\Delta_i\Delta_j G| )\right) &\leq& M.
\end{eqnarray*}
It is trivial to see that,
\begin{eqnarray}\label{convex1}
2^{-q} \sum_{j=1}^{2^q}\Psi_1\left(\sum_{i=1}^{2^p} \Phi_1 (|\Delta_i\Delta_j G| )\right) &\leq& 2^{-q} M.
\end{eqnarray}
As $\Psi_1$ is convex, so
\begin{eqnarray}\label{convex2}
\hskip -0.5cm2^{-q} \sum_{j=1}^{2^q}\Psi_1\left(\sum_{i=1}^{2^p} \Phi_1 (|\Delta_i\Delta_j G| )\right) &\geq& \Psi_1 \left(2^{-q} \sum_{j=1}^{2^q} \sum_{i=1}^{2^p} \Phi_1 (|\Delta_i\Delta_j G| )\right).
\end{eqnarray}
It turns out from (\ref{convex1}) and (\ref{convex2}) that
\begin{eqnarray}
\Psi_1 \left(2^{-q} \sum_{j=1}^{2^q} \sum_{i=1}^{2^p} \Phi_1 (|\Delta_i\Delta_jG| )\right)&\leq& 2^{-q}M.
\end{eqnarray}
This leads to
\begin{eqnarray}\label{n2}
    2^{-q} \sum_{j=1}^{2^q} \sum_{i=1}^{2^p} \Phi_1 (|\Delta_i\Delta_j G| )&\leq& \psi_1(2^{-q} M).
    \end{eqnarray}
    This is equivalent to
    \begin{eqnarray}\label{convex3}
2^{-p}2^{-q} \sum_{j=1}^{2^q} \sum_{i=1}^{2^p} \Phi_1 (|\Delta_i\Delta_j G| )&\leq& 2^{-p}\psi_1(2^{-q} M).
    \end{eqnarray}
    But, by the convexity of $\Phi_1$, we have
    \begin{eqnarray}\label{convex4}
    2^{-p}2^{-q} \sum_{j=1}^{2^q} \sum_{i=1}^{2^p} \Phi_1 (|\Delta_i\Delta_j G| )&\geq& \Phi_1 \left(2^{-p}2^{-q} \sum_{j=1}^{2^q} \sum_{i=1}^{2^p}  |\Delta_i\Delta_jG|
    \right).
    \end{eqnarray}
    So it follows from (\ref{convex3}) and (\ref{convex4}) that
    \begin{eqnarray}
    \Phi_1\left(2^{-p}2^{-q} \sum_{j=1}^{2^q} \sum_{i=1}^{2^p}  |\Delta_i\Delta_j G| \right)&\leq& 2^{-p}\psi_1 (2^{-q}M).
    \end{eqnarray}
    Therefore,
    \begin{eqnarray}
    2^{-p}2^{-q} \sum_{j=1}^{2^q} \sum_{i=1}^{2^p}  |\Delta_i\Delta_j G| &\leq& \varphi_1(2^{-p}\psi_1(2^{-q} M)).
    \end{eqnarray}
          So
    \begin{eqnarray}
    \sum_{j=1}^{2^q} \sum_{i=1}^{2^p}  |\Delta_i\Delta_j G| &\leq& 2^{p+q} \varphi_1(2^{-p}\psi_1(2^{-q} M)).
    \end{eqnarray}
    By the same method, one can see that
    \begin{eqnarray}\label{equ3}
\hskip -0.5cm\sum_{j=1,3,5,\cdots 2^{q+1}-1} \sum_{i=1,3,5,\cdots,2^{p+1}-1}  |\Delta_{i+1}\Delta_{j+1}G| \leq  2^{p+q}\varphi_1(2^{-p}\psi_1(2^{-q} M)).
    \end{eqnarray}
Therefore, it follows from (\ref{equ1}), (\ref{equ2}) and (\ref{equ3}) that there exists $K>0$ such that
\begin{eqnarray*}
&&\left|S(E_{p+1},E'_{q+1})-S(E_{p+1},E_q^{'})-S(E_{p},E_{q+1}^{'})+S(E_{p},E_q^{'})\right|\\
&\leq & K2^{p+q}\varphi_1(2^{-p}\psi_1(2^{-q}M))\varrho[\varphi(P 2^{-p} )]\sigma[\psi (Q 2^{-q})].
\end{eqnarray*}
\vskip 5pt
{\bf Step 2}: Let's prove that
\begin{eqnarray}
\lim_{p,q\to \infty} \left [S(E+E_p, E'+E'_q)-S(E_p, E'_q)\right ]=0.
\end{eqnarray}
Denoting $x_l$, $l=0,1,\cdots, L$ ($y_n$, $n=0,1,\cdots,L'$) the
distinct points of $E_p$ ($E'_q$) in increasing order, and by
$x_{l-1, i},\ i=0,1,\cdots, M_l$ ($y_{n-1, j},\ j=0,1,\cdots,
M'_n$) those of $E+E_p$ ($E'+E'_q$) lying in the interval
$x_{l-1}\leq x\leq x_l$ ($y_{n-1}\leq y\leq y_n$) with
$x_{l-1,0}=x_{l-1}, x_{l-1,M_l}=x_l$ ($y_{n-1,0}=y_{n-1},
y_{n-1,M'_n}=y_n$), so
\begin{eqnarray*}
&&|S(E+E_p,E'+E'_q)-S(E_p, E'_q)|\\
&=&|S(E+E_p,E'+E'_q)-S(E+E_p,E'_q))+(S(E+E_p,E'_q)-S(E_p, E'_q)|\\
&=&\left |\sum_{l=1}^L\sum_{i=1}^{M_l}\sum_{n=1}^{L'}\sum_{j=1}^{M'_n}\bigg\{[F(x_{l-1,i-1},y_{n-1,j-1})-F(x_{l-1,i-1},y_{n-1})]\right .\\
&&\hskip 1cm+[F(x_{l-1,i-1},y_{n-1})-F(x_{l-1},y_{n-1})]\bigg\}\bigg\{G(x_{l-1,i},y_{n-1,j})\\
&&\left .\hskip 1cm-G(x_{l-1,i-1},y_{n-1,j})
-G(x_{l-1,i},y_{n-1,j-1})+G(x_{l-1,i-1},y_{n-1,j-1})\bigg\}\right |\\
&\leq & 4N_1N_2[\psi (2^{-q}Q)+\varphi(2^{-p}P)]\cdot \max |G|\\
&\to & 0 ,\  as\ p,\ q\to \infty.
\end{eqnarray*}
Here $N_1,N_2$ denote the number of points of $E+E_0$, $E'+E'_0$, respectively.
\vskip 5pt
{\bf Step 3}: Let $F(x,y)$ vanish for $x=x'$ identically in $y$, and for $y=y'$ identically in $x$, so
\begin{eqnarray*}
&&F_{E_0, E'}(x,y)=F(x',y)=0,\ \ F_{E, E'_0}(x,y)=F(x,y')=0,\\
&&F_{E_0, E'_0}(x,y)=F(x',y')=0.
\end{eqnarray*}
If this is so, and also note that $S(E,E')=S_{F_{E,E'}}(E+E_p,E'+E'_q)$, then from Step 2, Step 1
and (\ref{hz38}),
\begin{eqnarray}\label{new2}
&&
\left|S(E,E')\right|\nonumber\\
&=&\left|S(E,E')-S(E_0,E')-S(E,E'_0)+S(E_0,E'_0)\right|\nonumber\\
&=&\Big|\lim_{p,q\to \infty}\Big[S_{F_{E,E'}}(E+E_p, E'+E'_q)-S_{F_{E,E'}}(E_p, E'_q)+S_{F_{E,E'}}(E_p,E'_q)\nonumber\\
&&\hskip 1cm-S_{F_{E,E'}}(E_0, E_q')-S_{F_{E,E'}}(E_p,E'_0)+S_{F_{E,E'}}(E_0,E'_0)\Big]\Big|\nonumber\\
&=&\Big|\lim_{p,q\to \infty}\Big[S_{F_{E,E'}}(E_p,E'_q)-S_{F_{E,E'}}(E_0, E_q')-S_{F_{E,E'}}(E_p,E'_0)\nonumber\\
&&\hskip 2cm+S_{F_{E,E'}}(E_0,E'_0)\Big]\Big|\nonumber\\
&=&\Big|\sum_{p,q=0}^{\infty}\Big[S_{F_{E,E'}}(E_{p+1},E'_{q+1})-S_{F_{E,E'}}(E_{p+1},E'_q)
\nonumber\\
&&\hskip2cm -S_{F_{E,E'}}(E_{p},E'_{q+1})+S_{F_{E,E'}}(E_{p},E'_{q})\Big]\Big|\nonumber\\
&\leq&\sum_{p,q=0}^{\infty}K2^{p+q}\varrho[\varphi(P 2^{-p} )]\sigma[\psi (Q 2^{-q})]\varphi_1(2^{-p}\psi_1(2^{-q}M))\nonumber\\
&\leq&4K\sum_{m,n=1}^{\infty}\varrho[\varphi ({P\over n})]\sigma[\psi ({Q\over m})]\varphi_1({1\over n}\psi_1({1\over m} M)).
\end{eqnarray}
Let $F_{x',y'}(x,y):= F(x,y)-F(x',y)-F(x,y')+F(x',y')$ and replace
$F(x,y)$ by $F_{x',y'}(x,y)$ for $x'\leq x\leq x''$, $y'\leq y
\leq y''$. This alteration doesn't affect double difference of
$F$. Therefore we may suppose that $F(x,y)$ vanish identically on
the lines $x=x'$ and $y=y'$ as above. \vskip 5pt

{\bf Step 4}: We determine a set of finite points $H_v\times
H_{v'}:=\{x'=x_0<x_1<\cdots<x_L=x'',y'=y_0<y_1<\cdots<y_{L'}=y''\}$,
where $L\leq 2\cdot2^v$, $L'\leq 2\cdot2^{v'}$, such that in the
rectangle
$[x_{l-1}+\delta,x_l-\delta]\times[y_{k-1}+\delta,y_k-\delta]$,
$|\Delta_i\Delta_j G|<\varepsilon(v,v')$ for any $0<\delta \leq
{1\over 4} \min\{\min_{1\leq l\leq L}\{x_l-x_{l-1}\},  \min_{1\leq
l\leq L^{\prime}}\{y_l-y_{l-1}\}\}$. Moreover, in the interval
$[x',x'']\times[y_{k-1}+\delta,y_k-\delta]$, the total
$\Psi$-variation of $F$ in $y$ is at most $Q\cdot 2^{-v'}$, the
total $\Phi_1$, $\Psi_1-$ variation of $G$ is at most $M\cdot
2^{-v'}$; in the interval
$[x_{l-1}+\delta,x_l-\delta]\times[y',y'']$, the total
$\Phi$-variation of $F$ in $x$ is at most $P\cdot 2^{-v}$ and the
total $\Phi_1$-variation of $G$ in $x$ for the given partition of
$H_v^{\prime}$ of $[y',y'']$ is at most $2^{-v}L'\psi_1({1\over
L'}M)$. Here, the first and the third statements are obvious, the
second statement follows from the finite large jump condition. The
last one
 can be seen by observing that
$\sum_{j=1}^{L'}\Psi_1\left(\sum_{i=1}^L\Phi_1(|\Delta_i\Delta_j
G|)\right)\leq M$ is equivalent to
$\sum_{j=1}^{L'}\sum_{i=1}^L\Phi_1(|\Delta_i\Delta_j G|)\leq
L'\psi_1({1\over L'}M)$. More generally, for any partition
$E=\{x_1^{\prime},x_2^{\prime},\cdots,x_N^{\prime}\}$, we have
$\sum_{i=1}^{N} \sum_{j=1}^{L'}\Phi_1(|\Delta_i^{\prime}\Delta_j
G|)\leq L'\psi_1({1\over L'}M)$. Here $\Delta_i^{\prime}\Delta_j
G$ is the double increment of $G$ on
$(x_{i-1}^{\prime},x_i^{\prime})\times (y_{j-1},y_j)$. We can make
$E$ include $H_v$ among their points of divisions and denote
$E_{l-1}=\{x_{l-1,1}, x_{l-1,2}, \cdots , x_{l-1,N_{l-1}}\}$ all
the points in $E$ falling into the interval $(x_{l-1},x_l)$
$(l=1,2,\cdots, L)$. We can certainly make
\begin{eqnarray}\label{hz34}
\sum_{i=1}^{N_{l-1}}
\sum_{j=1}^{L'}\Phi_1(|\Delta_{l-1,i}^{\prime}\Delta_j G|)\leq
2^{-v}L'\psi_1({1\over L'}M),
\end{eqnarray}
where $\Delta_{l-1,i}^{\prime}\Delta_j G$ is the double increment of $G$
on $(x_{l-1,i-1}^{\prime},x_{l-1,i}^{\prime})\times (y_{j-1},y_j)$. In fact $E_{l-1}$ can be any partition of
$[x_{l-1}+\delta,x_l-\delta]$ for any sufficiently small $\delta >0$.

We need to prove that for any $\varepsilon>0$,
\begin{eqnarray}\label{n1}
|S(\tilde{D},\tilde{D'})-S(D,D')|<\varepsilon,
\end{eqnarray}
as long as $\tilde{D}\times\tilde{D'}$, and $D\times D'$ include
$H_v\times H_{v'}$. Observe that
\begin{eqnarray}\label{hz35}
&&|S(\tilde{D},\tilde{D'})-S(D,D')|\nonumber\\
&\leq&
|S(\tilde{D},\tilde{D'})-S(\tilde{D},D')|+|S(\tilde{D},D')-S(D,D')|\nonumber\\
&\leq&|\int_{x'}^{x''}\int_{y'}^{y''}(F_{\tilde{D},\tilde{D'}}-F_{\tilde
D,D'})d_{x,y}G(x,y)|\nonumber\\
&&+|\int_{x'}^{x''}\int_{y'}^{y''}(F_{\tilde{D},{D'}}-F_{
D,D'})d_{x,y}G(x,y)|.
\end{eqnarray}
First, since $F_{\tilde{D},\tilde{D'}}-F_{\tilde D,D'}$ vanishes
identically in $x$, when $y=y_{k-1}$, from Step 3 and (\ref{hz38}), (\ref{hz39}),  we obtain
for any sufficiently small $\delta>0$,
\begin{eqnarray}\label{hz36}
&&|S(\tilde{D},\tilde{D'})-S(\tilde{D},D')|\nonumber\\
&\leq&
\sum_{k=1}^{L'}\left|\int_{x'}^{x''}\int_{y_{k-1}+\delta}^{y_k-\delta}
(F_{\tilde{D},\tilde{D'}}-F_{\tilde D,D'})d_{x,y}G(x,y)\right|\nonumber\\
&\leq& \sum_{k=1}^{L'}4K\sum_{m,n=1}^{\infty}\varrho[\varphi
({P\over n})]\sigma[\psi ({2^{-v'}Q\over m})]\varphi_1[{1\over
n}\psi_1({1\over m} 2^{-v'}M)]\nonumber\\
&\leq&
\sum_{k=1}^{L'}4K\sum_{n=1}^{\infty}\sum_{q=0}^{\infty}\varrho[\varphi
({P\over n})]2^q\sigma[\psi ({2^{-(v'+q)}Q})]\varphi_1[{1\over n}\psi_1( 2^{-(v'+q)}M)]\nonumber\\
&=& L'4K\sum_{n=1}^{\infty}\sum_{q={v'}}^{\infty}\varrho[\varphi
({P\over n})]2^{-{v'}+q}\sigma[\psi ({2^{-q}Q})]\varphi_1[{1\over n}\psi_1( 2^{-q}M)]\nonumber\\
&\leq&16K\sum_{n=1}^{\infty}\sum_{m=2^{v'-1}+1}^{\infty}\varrho[\varphi
({P\over n})]\sigma[\psi ({Q\over m})]\varphi_1[{1\over
n}\psi_1({1\over m}M)]\nonumber\\
&\leq&\varepsilon_{v'},
\end{eqnarray}
where $\varepsilon_{v'}\to 0$, as ${v'}\to \infty$.

Second, since $F_{\tilde D,D'}-F_{D,D'}$ vanishes identically in $y$,
when $x=x_{l-1}$. From the discussion above and (\ref{hz34}), we know for any partition $E_p=\{x_1,x_2,
\cdots, x_{2^{p}}\}$ of $[x_{l-1}+\delta,x_l-\delta]$, and any partition $E_q{\prime}=\{y_1,y_2,
\cdots, y_{2^{q}}\}$ of $[y^{\prime},y^{\prime\prime}]$,
 (\ref{n2}) bacomes
 \begin{eqnarray}\label{n3}
    2^{-q} \sum_{j=1}^{2^q} \sum_{i=1}^{2^p} \Phi_1 (|\Delta_i\Delta_j G| )&\leq& 2^{-v}\psi_1(2^{-q}
    M).
    \end{eqnarray}
So from Step 3 and (\ref{hz38}), (\ref{hz39}),
\begin{eqnarray*}
&&\left|\int_{x_{l-1}+\delta}^{x_l-\delta}\int_{y'}^{y''}
(F_{\tilde{D},D'}-F_{D,D'})d_{x,y}G(x,y)\right|\\
&\leq& 4K\sum_{m,n=1}^{\infty}\varrho[\varphi
({2^{-v}P\over n})]\sigma[\psi ({Q\over
m})]\varphi_1[{{2^{-v}}\over
n}\psi_1({1\over m}M)]\\
&\leq & 4K\sum_{m=1}^{\infty}\sum_ {q=0}^{\infty}2^q\varrho[\varphi
({2^{-(v+q)}P})]\sigma[\psi ({Q\over
m})]\varphi_1[{{2^{-(v+q)}}}\psi_1({1\over m}M)]\\
&= & 4K2^{-v}\sum_{m=1}^{\infty}\sum_ {q=v}^{\infty}2^q\varrho[\varphi
({2^{-q}P})]\sigma[\psi ({Q\over
m})]\varphi_1[{2^{-q}}\psi_1({1\over m}M)]\\
&\leq& 8K2^{-v}\sum_{m=1}^{\infty}\sum_{n=2^{v-1}+1}^{\infty}\varrho[\varphi
({P\over n})]\sigma[\psi ({Q\over m})]\varphi_1[{1\over
n}\psi_1({1\over m}M)].
\end{eqnarray*}
Now it turns out that
\begin{eqnarray}\label{hz37}
&&|S(\tilde{D},D')-S(D,D')|\nonumber\\
&\leq&
\sum_{l=1}^{L}\left|\int_{x_{l-1}+\delta}^{x_l-\delta}\int_{y'}^{y''}
(F_{\tilde{D},D'}-F_{D,D'})d_{x,y}G(x,y)\right|\nonumber\\
&\leq&16K\sum_{m=1}^{\infty}\sum_{n=2^{v-1}+1}^{\infty}\varrho[\varphi
({P\over n})]\sigma[\psi ({Q\over m})]\varphi_1[{1\over
n}\psi_1({1\over m}M)]\nonumber\\
&\leq&\varepsilon_{v},
\end{eqnarray}
where $\varepsilon_{v}\to 0$, as ${v}\to \infty$.

Thus we can get (\ref{n1}) from (\ref{hz35}), (\ref{hz36}) and (\ref{hz37}), as $v,v'\to \infty$, which means
$S(D,D')$ is a Cauchy sequence, so $\lim\limits_{m(D\times D')\to
0}S(D,D')$ exists. In the following, we show the limit is unique.
For this, let $D_1\times D'_1$, $D_2\times D'_2$ be arbitrary two
partitions of $[x',x'']\times[y',y'']$ including $H\times H'$.
From the above we know,
\begin{eqnarray*}
&&|S(D_1\cup D_2,D'_1\cup D'_2)-S(D_1,D'_1)|\to 0, \ as \  m(D_1\times D'_1)\to 0,\\
&&|S(D_1\cup D_2,D'_1\cup D'_2)-S(D_2,D'_2)|\to 0, \ as \
m(D_2\times D'_2)\to 0.
\end{eqnarray*}
Therefore,
\begin{eqnarray*}
&&\lim_{m(D_1\times D'_1)\to 0}S(D_1,D'_1)=\lim_{m(D_2\times D'_2)\to 0}S(D_2,D'_2)\\
&&=\lim_{m(D_1\times D'_1),m(D_2\times D'_2)\to 0}S(D_1\cup
D_2,D'_1\cup D'_2),
\end{eqnarray*}
that is to say, $\lim\limits_{m(D\times D')\to 0}S(D,D')$ is
unique, and we define it as\\
$\int_{x'}^{x''}\int_{y'}^{y''}F(x,y)d_{x,y} G(x,y)$. So we proved
our theorem. $\hfill\diamond$
\bigskip

 In the following we say an integral is well defined which is
in the sense of Theorem \ref{ffcr1}.
 The following convergence theorem plays an important role in
establishing It${\hat {\rm o}}$'s formula:
\begin{thm}\label{our3}
Assume there exist convex functions $\Phi$, $\Psi$, $\Phi_1$,
$\Psi_1$ such that $F_k(x,y)$ and $F(x,y)$ are continuous
functions of bounded $\Phi$-variation in $x$ uniformly in $y$ and
of bounded $\Psi$-variation in $y$ uniformly in $x$ and for $F_k$
uniformly in $k$, $G(x,y)$ and $G_k(x,y)$ are of bounded $\Phi_1$,
$\Psi_1-$variation in ($x,y$) uniformly in $k$ and satisfies the finite large
jump condition. If there exist
$\varrho_i$ and $\sigma_i$ subject to $\varrho_i(u)\sigma_i(u)=u$,
i=1,2, and a small positive number $\delta>0$ such that
\begin{eqnarray}\label{assumption1}
&&\sum\limits_{m,n}\varrho_1[\varphi({1\over n})]\sigma_1[\psi({1\over m})]\varphi_1^{1\over {1+\delta}}[{1\over n}\psi_1({1\over m})]\nonumber\\
&&+
\sum\limits_{m,n}\varrho_2[\varphi(({1\over n})^{1\over {1+\delta}})]\sigma_2[\psi({1\over m})]\varphi_1[{1\over n}\psi_1({1\over m})]<\infty,
\end{eqnarray}
or
\begin{eqnarray}\label{assumption2}
&&\sum\limits_{m,n}\varrho_1[\varphi({1\over n})]\sigma_1[\psi({1\over m})]\varphi_1^{1\over {1+\delta}}[{1\over n}\psi_1({1\over m})]\nonumber\\
&&+
\sum\limits_{m,n}\varrho_2[\varphi({1\over n})]\sigma_2[\psi(({1\over m})^{1\over {1+\delta}})]\varphi_1[{1\over n}\psi_1({1\over m})]<\infty,
\end{eqnarray}
 and let $F_k(x,y)\to F(x,y)$, $G_k(x,y)\to G(x,y)$ as $k\to \infty$ uniformly in ($x,y$). Then we have
\begin{eqnarray}
\int_{x'}^{x''}\int_{y'}^{y''}F_k(x,y)d_{x,y}G_k(x,y)\to\int_{x'}^{x''}\int_{y'}^{y''}F(x,y)d_{x,y}G(x,y),
\end{eqnarray}
when $k\to \infty$.
\end{thm}
{\bf Proof:} First note that from Theorem \ref{ffcr1} under the above assumptions,
the integral $\int_{x'}^{x''}\int_{y'}^{y''}F_k(x,y)d_{x,y}G_k(x,y)$ and
$\int_{x'}^{x''}\int_{y'}^{y''}F(x,y)d_{x,y}G(x,y)$ are all well defined. It's easy to see
that
\begin{eqnarray*}
&&{1 \over 2}\left(\int_{x'}^{x''}\int_{y'}^{y''}F_k(x,y)d_{x,y}G_k(x,y)-\int_{x'}^{x''}\int_{y'}^{y''}F(x,y)d_{x,y}G(x,y)\right)\\
&=& \int_{x'}^{x''}\int_{y'}^{y''}F_k(x,y)d_{x,y}{1\over 2}(G_k(x,y)-G(x,y))\\
&&+\int_{x'}^{x''}\int_{y'}^{y''}{1\over
2}(F_k(x,y)-F(x,y))d_{x,y}G(x,y).
\end{eqnarray*}
We study ${1\over 2}$ of the integral only for convenience in what follows.
First consider the integral $\int_{x'}^{x''}\int_{y'}^{y''}F_k(x,y)d_{x,y}(G_k(x,y)-G(x,y))$.
Note there exist constant $P_1, Q_1, M_1, M_2>0$, which are independent of $k$ such that for
any partition $E\times E'$defined before
\begin{eqnarray}
&&\sum_{i=1}^N\Phi(|F_k(x_i,y)-F_k(x_{i-1},y)|)\leq P_1,\ \ \   for \ any\ y\in[y',y''],\\
&&\sum_{j=1}^{N'}\Psi(|F_k(x,y_j)-F_k(x,y_{j-1})|)\leq Q_1,\ \ \   for \ any\ x\in[x',x''],\\
&&\sum_{j=1}^{N'}\Psi_1\left(\sum_{i=1}^N\Phi_1(|\Delta_i\Delta_j G_k|)\right)\leq M_1,\\
&&\sum_{j=1}^{N'}\Psi_1\left(\sum_{i=1}^N\Phi_1(|\Delta_i\Delta_j G|)\right)\leq M_2.
\end{eqnarray}
For the small $\delta>0$ given in condition (\ref{assumption1}),
from the convexity of $\Phi_1$ and $\Psi_1$ and $G_k\to G$ when
$k\to \infty$, we have
\begin{eqnarray*}
&&\sum_{j=1}^{N'}\Psi_1\left(\sum_{i=1}^N\Phi_1\Big(|\Delta_i\Delta_j {1\over 2}(G_k-G)|^{1+\delta}\Big)\right)\\
&=&\sum_{j=1}^{N'}\Psi_1\left(\sum_{i=1}^N\Phi_1\Big(|\Delta_i\Delta_j {1\over 2}(G_k-G)|^{\delta}\cdot|\Delta_i\Delta_j {1\over 2}(G_k-G)|\Big)\right)\\
&\leq& \sum_{j=1}^{N'}\Psi_1\left(\sum_{i=1}^N|\Delta_i\Delta_j {1\over 2}(G_k-G)|^{\delta}\Phi_1\Big(|\Delta_i\Delta_j {1\over 2}(G_k-G)|\Big)\right)\\
&\leq &\sum_{j=1}^{N'}\Psi_1\left(\max_i|\Delta_i\Delta_j {1\over 2}(G_k-G)|^{\delta}\sum_{i=1}^N\Phi_1\Big(|\Delta_i\Delta_j {1\over 2}(G_k-G)|\Big)\right)\\
&\leq &\sum_{j=1}^{N'}\max_i|\Delta_i\Delta_j {1\over 2}(G_k-G)|^{\delta}\Psi_1\left(\sum_{i=1}^N\Phi_1\Big({1\over 2}|\Delta_i\Delta_j G_k|+{1\over 2}|\Delta_i\Delta_j G|\Big)\right)\\
&\leq &\max_{i,j}|\Delta_i\Delta_j (G_k-G)|^{\delta}\sum_{j=1}^{N'}\Psi_1\left(\sum_{i=1}^N\left ({1\over 2}\Phi_1\Big(|\Delta_i\Delta_j G_k|\Big)+{1\over 2}\Phi_1\Big(|\Delta_i\Delta_j G|\Big)\right)\right )\\
&\leq &\varepsilon_1(k) M,
\end{eqnarray*}
where $\varepsilon_1(k)\to 0$ as $k\to \infty$, and $M$ is a constant independent of $k$. If we define
\begin{eqnarray*}
S(E,E')=\sum_{i=1}^N\sum_{j=1}^{N'} F_k(x_{i-1},y_{j-1})\Big(\Delta_i\Delta_j (G_k-G)\Big),
\end{eqnarray*}
and similar to (\ref{new2}), by dominated convergence theorem to
the infinite series,
\begin{eqnarray*}
|S(E,E')|\leq 4K\sum\limits_{m,n}\varrho_1[\varphi({P_1\over
n})]\sigma_1[\psi({Q_1\over m})]\varphi_1^{1\over
{1+\delta}}[{1\over n}\psi_1({2\varepsilon_1(k)\over m}M)] \to 0,\
as\ k\to \infty,
\end{eqnarray*}
as the series $\sum\limits_{m,n}\varrho_1[\varphi({1\over n})]\sigma_1[\psi({1\over m})]\varphi_1^{1\over {1+\delta}}[{1\over n}\psi_1({1\over m})]< \infty$.
This implies as $k\to \infty$,
\begin{eqnarray}\label{limit1}
\lim_{k\to \infty}\int_{x'}^{x''}\int_{y'}^{y''}F_k(x,y)d_{x,y}(G_k(x,y)-G(x,y))=0.
\end{eqnarray}
For the second integral $\int_{x'}^{x''}\int_{y'}^{y''}(F_k(x,y)-F(x,y))d_{x,y}G(x,y)$, we can use a similar method to prove
\begin{eqnarray}\label{limit2}
\lim\limits_{k\to\infty}\int_{x'}^{x''}\int_{y'}^{y''}(F_k(x,y)-F(x,y))d_{x,y}G(x,y)=0.
\end{eqnarray}
For this, we note from the assumption there is a $\delta >0$ such that,
\begin{eqnarray*}
&&\sum_{i=1}^N \Phi^{1+\delta}\Big(|{1\over 2}(F_k-F)(x_i,y)-{1\over 2}(F_k-F)(x_{i-1},y)|\Big)\\
&\leq &\max_i \Phi^{\delta}\Big(|{1\over 2}(F_k-F)(x_i,y)-{1\over 2}(F_k-F)(x_{i-1},y)|\Big)\cdot\\
&&\sum_{i=1}^N {1\over 2}\Phi\Big((|F_k(x_i,y)-F_k(x_{i-1},y)|)+\Phi(|F_k(x_i,y)-F_k(x_{i-1},y)|)\Big) \\
&\leq &\varepsilon_2(k)M',
\end{eqnarray*}
where $\varepsilon_2(k)\to 0$, as $k\to \infty$, and $M'$ is a constant independent of $k$.
So under the assumption $ \sum\limits_{m,n}\varrho_2[\varphi(({1\over n})^{1\over {1+\delta}})]\sigma_2[\psi({1\over m})]\varphi_1[{1\over n}\psi_1({1\over m})]<\infty$, we can prove (\ref {limit2}) using the same argument in proving (\ref{limit1}). Therefore under assumption (\ref{assumption1}), we prove the desired result. The proof is similar under the assumption  (\ref{assumption2}).
$\hfill\diamond$
\vskip5pt
\begin{rmk}
From the proof we can easily see that under the condition that there exist two functions $\varrho$ and $\sigma$ subject to $\varrho(u)\sigma(u)=u$ and a small number $\delta>0$ such that
\begin{eqnarray}
&&\sum\limits_{m,n}\varrho[\varphi({1\over n})]\sigma[\psi({1\over m})]\varphi_1^{1\over {1+\delta}}[{1\over n}\psi_1({1\over m})]<\infty.
\end{eqnarray}
Then as $k\to \infty$,
\begin{eqnarray}
\int_{x'}^{x''}\int_{y'}^{y''}F(x,y)d_{x,y}G_k(x,y)\to\int_{x'}^{x''}\int_{y'}^{y''}F(x,y)d_{x,y}G(x,y).
\end{eqnarray}
Similarly, under the condition that there exist two functions $\varrho$ and $\sigma$ subject to $\varrho(u)\sigma(u)=u$ and a small number $\delta>0$ such that
\begin{eqnarray}
\sum\limits_{m,n}\varrho[\varphi(({1\over n})^{1\over {1+\delta}})]\sigma[\psi({1\over m})]\varphi_1[{1\over n}\psi_1({1\over m})]<\infty,
\end{eqnarray}
or
\begin{eqnarray}
\sum\limits_{m,n}\varrho[\varphi({1\over n})]\sigma[\psi(({1\over m})^{1\over {1+\delta}})]\varphi_1[{1\over n}\psi_1({1\over m})]<\infty.
\end{eqnarray}
Then as $k\to \infty$,
\begin{eqnarray}
\int_{x'}^{x''}\int_{y'}^{y''}F_k(x,y)d_{x,y}G(x,y)\to\int_{x'}^{x''}\int_{y'}^{y''}F(x,y)d_{x,y}G(x,y).
\end{eqnarray}
\end{rmk}

It is easy to see that in the definition of
$\int_{x'}^{x''}\int_{y'}^{y''}F(x,y)d_{x,y}G(x,y)$, one can take
$F(x_i,y_j)$ instead of $F(x_{i-1},y_{j-1})$ in (\ref{chunrong1}).
One can also prove the convergence of (\ref{chunrong1}) in this
case and denote the integral by
$\int_{x'}^{x''}\int_{y'}^{y''}F(x,y)d_{x,y}^*G(x,y)$, the
backward integral. In general, this should be different from
$\int_{x'}^{x''}\int_{y'}^{y''}F(x,y)d_{x,y}G(x,y)$. But under
slightly stronger conditions than those in Theorem \ref{ffcr1}, as
in the one-parameter case, these two integrals equal. This result
is proved in the following proposition.

\begin{prop}\label{proposition}
Assume there exist convex functions $\Phi$, $\Psi$, $\Phi_1$,
$\Psi_1$ such that $F(x,y)$ is a continuous function of bounded
$\Phi$-variation in $x$ uniformly in $y$ and of bounded
$\Psi$-variation in $y$ uniformly in $x$,
 $G(x,y)$ is of
bounded $\Phi_1,\Psi_1$-variation in $(x,y)$ and satifies the finite
large jump condition. If there exist
functions $\varrho$ and $\sigma$ subject to
$\varrho(u)\sigma(u)=u$ and a small positive $\delta>0$ such that
one of the following two conditions is satisfied

(i) 
\begin{eqnarray*}
\sum\limits_{m,n}\varrho[\varphi(({1\over n})^{1\over
{1+\delta}})]\sigma[\psi({1\over m})]\varphi_1[{1\over
n}\psi_1({1\over m})]<\infty,
\end{eqnarray*}

(ii)
\begin{eqnarray*}
\sum\limits_{m,n}\varrho[\varphi({1\over n})]\sigma[\psi(({1\over
m})^{1\over {1+\delta}})]\varphi[{1\over n}\psi_1({1\over
m})]<\infty.
\end{eqnarray*}
Then
\begin{eqnarray*}
\int_{x'}^{x''}\int_{y'}^{y''}F(x,y)d_{x,y}G(x,y)=\int_{x'}^{x''}\int_{y'}^{y''}F(x,y)d^*_{x,y}G(x,y).
\end{eqnarray*}
\end{prop}
{\bf Proof}: We only prove the result when condition (i) is satisfied. Denote
\begin{eqnarray*}
&&S(E,E')=\sum_{i=1}^N\sum_{j=1}^{N'}F(x_{i-1},y_{j-1})\Delta_i\Delta_j
G,\\
&&S^*(E,E')=\sum_{i=1}^N\sum_{j=1}^{N'}F(x_{i},y_{j})\Delta_i\Delta_j
G.
\end{eqnarray*}
Here $E$ and $E'$ are the same as before. Denote
\begin{eqnarray*}
\tilde{F}_{\delta_{x_{i-1}},\delta_{y_{j-1}}}(x_{i-1},y_{j-1})
=F(x_{i-1}+\delta_{x_{i-1}},y_{j-1}+\delta_{y_{j-1}})-F(x_{i-1},y_{j-1}).
\end{eqnarray*}
Here $\delta_{x_{i-1}}=x_i-x_{i-1}$,
$\delta_{y_{j-1}}=y_j-y_{j-1}$. Then
\begin{eqnarray*}
S^*(E,E')-S(E,E')=2\sum_{i=1}^N\sum_{j=1}^{N'}{1\over
2}\tilde{F}_{\delta_{x_{i-1}},\delta_{y_{j-1}}}(x_{i-1},y_{j-1})
\Delta_i\Delta_j G.
\end{eqnarray*}
Note from the assumptions, there is a $\delta>0$ such that
\begin{eqnarray*}
&&\sum_{i=1}^N \Phi^{1+\delta}\bigg(\Big|{1\over
2}[\tilde{F}_{\delta_{x_{i}},\delta_{y_{j-1}}}(x_{i},y_{j-1})-\tilde{F}_{\delta_{x_{i-1}},\delta_{y_{j-1}}}(x_{i-1},y_{j-1})]\Big|\bigg)\\
 &\leq& \max_{i,j}\Phi^{\delta}\bigg({1\over
2}\Big|\big[F(x_{i+1},y_{j})-F(x_{i},y_{j-1})\big]
 -\big[F(x_{i},y_{j})-F(x_{i-1},y_{j-1})\big]\Big|\bigg)\\
 &&\sum_{i=1}^N
\Phi\bigg(\Big|{1\over 2}\big[F(x_{i+1},y_{j})-F(x_{i},y_{j})\big]
-{1\over
2}\big[F(x_{i},y_{j-1})-F(x_{i-1},y_{j-1})\big]\Big|\bigg)\\
&\leq& \max_{i,j}\Phi^{\delta}\bigg({1\over
2}\Big|\big[F(x_{i+1},y_{j})-F(x_{i},y_{j})\big]
 -\big[F(x_{i},y_{j-1})-F(x_{i-1},y_{j-1})\big]\Big|\bigg)\\
 &&\sum_{i=1}^N {1\over 2}\Big(\Phi\big(|F(x_{i+1},y_{j})-F(x_{i},y_{j})|\big)+\Phi\big(|F(x_{i},y_{j-1})-F(x_{i-1},y_{j-1})|\big)\Big)\\
&\leq& \varepsilon(E,E')P,
\end{eqnarray*}
where $\varepsilon(E,E')\to 0$, when $m(E,E')\to 0$ and $P$ is a
constant. Therefore following (\ref{new2}), we see that
\begin{eqnarray*}
&&|S^*(E,E')-S(E,E')|\\
&=&
|\sum_{i=1}^N\sum_{j=1}^{N'}\tilde{F}_{\delta_{x_{i-1}},\delta_{y_{j-1}}}(x_{i-1},y_{j-1})
\Delta_i\Delta_j G|\\
&\leq& 8K
\sum_{m,n=1}^{\infty}\varrho[\varphi\Big(\big({2\varepsilon(E,E')P\over
n}\big)^{1\over {1+\delta}}\Big)]\sigma[\psi({Q\over
m})]\varphi_1({1\over
n}\psi_1({1\over m}M))\\
&\to& 0, \ as\ \varepsilon(E,E')\to 0.
\end{eqnarray*}
Therefore
\begin{eqnarray*}
 S^*(E,E')-S(E,E')\to 0\ as\ \varepsilon(E,E')\to 0.
 \end{eqnarray*}
 That is to say,
\begin{eqnarray*}
\int_{x'}^{x''}\int_{y'}^{y''}F(x,y)d_{x,y}G(x,y)=\int_{x'}^{x''}\int_{y'}^{y''}F(x,y)d^*_{x,y}G(x,y).
\end{eqnarray*}
$\hfill\diamond$
\vskip5pt

From Theorem \ref{ffcr1} we can easily generalize it to the
multi-parameter integral.
\begin{defi}
Let $E_1\times \cdots \times E_n = \{ a_1=x_1^0<x_1^1<\cdots<x_1^{N_1}=b_1, \cdots, a_n=x_n^0<x_n^1<\cdots<x_n^{N_n}=b_n\} $
 be an arbitrary partition of $[a_1,b_1]\times\cdots[a_n,b_n]$. We call $F(x_1,\cdots, x_n)$ is of bounded $\Phi_i$-variation in $x_i$ uniformly in $(x_1,\cdots,x_{i-1},x_{i+1},\cdots, x_n)$, $i=1,\cdots n$, if
\begin{eqnarray}
\sup_{x_1,\cdots,x_{i-1},x_{i+1},\cdots, x_n} \sup_{E_i} \sum_{k_i=1}^{N_i}\Phi_i(|\Delta_{x_i^{k_i-1},x_i^{k_i}}F|)<\infty,
\end{eqnarray}
Here $\Delta$ is the difference operator (see \cite{ash}) as follows,
\begin{eqnarray*}
\Delta_{x_i^{k_i-1},x_i^{k_i}}F&=&F(x_1,\cdots,x_{i-1},x_i^{k_i},x_{i+1},\cdots,x_n)\nonumber\\
&&-F(x_1,\cdots,x_{i-1},x_i^{k_i-1},x_{i+1},\cdots,x_n).
\end{eqnarray*}
We call $G(x_1,\cdots,x_n)$ is of bounded $\Psi_1,\cdots,\Psi_n$-variation in $(x_1,\cdots,x_n)$, if
\begin{eqnarray}
\hskip -0.7cm\sup_{E_1\times\cdots\times
E_n}\sum_{k_n=1}^{N_n}\Psi_n\Big(\cdots\big(\sum_{k_1=1}^{N_1}\Psi_1(|\Delta_{x_n^{k_n-1},x_n^{k_n}}\cdots\Delta_{x_1^{k_1-1},x_1^{k_1}}G|)\big)\cdots\Big)<\infty.
\end{eqnarray}
\end{defi}

We say a function $f(x_1,\cdots,x_n)$  has a jump at $(x_1^0,\cdots, x_n^0)$ if there
exists an $\varepsilon >0$ such that for any $\delta>0$, there exists
$(x_1^1,\cdots, x_n^1)$ satisfying  $\max\{|x_1^0-x_1^1|,\cdots,|x_n^0-x_n^1|\}<\delta$
and
$|\Delta_{x_n^0,x_n^1},\cdots,\Delta_{x_1^0,x_1^1}f|>\varepsilon$. For a
function $G(x_1,\cdots,x_n)$ of bounded $\Psi_1,\cdots,
\Psi_n$-variation, for any
given $\varepsilon>0$, it is easy to see that there exists a $\delta (\varepsilon)>0$ and
 a finite number of jump points $\{(x_1^1,\cdots,x_n^1),\cdots,(x_1^{m_1},\cdots,x_n^{m_n})\}$ such that
  $|\Delta_{\tilde x_n,x_n},\cdots,\Delta_{\tilde x_1,x_1} G|<\varepsilon$
whenever
$\max\{|\tilde x_1-x_1|,\cdots ,|\tilde x_n-x_n|\}<\delta(\varepsilon)$,
 $[\tilde x_i,x_i]\cap\{x_1^i,\cdots,x_n^{m_i}\}= \emptyset$ for all $i=1,2,\cdots, n$.
Denote $H_{10}\times\cdots \times
H_{n0}:=\{x_1^1,\cdots,x_1^{m_1}\}\times\cdots \times\{x_n^{1},\cdots,x_n^{m_n}\}$.
In the following, we assume the following finite large jump
condition: for any $\varepsilon >0$, there exists at
most finite many points $ \{x_1^1,\cdots,x_1^{m_1'}\}, \cdots,\{x_n^{1},\cdots,x_n^{m_n'}\}$
and a constant $\delta (\varepsilon)
>0$ such that for each
$i=1,2,\cdots, n$, the total $\Psi_1,\cdots ,
\Psi_n$-variation of $G$ on $[x_1', x_1'']\times \cdots [x_i,x_i+\delta ]\times \cdots\times [x_n',x_n'']$ is smaller than $\varepsilon$
if $[x_i, x_i+\delta]\cap \{x_i^1, \cdots ,x_i^{m_i^\prime}\}= \emptyset$.
Denote $H_1\times
\cdots \times H_n:=\{x_1^1,\cdots,x_1^{m_1'}\}\times \cdots \times \{x_n^{1},\cdots,x_n^{m_n'}\}$. It is
obvious that $H_1\times\cdots \times H_n\supset H_{10}\times\cdots\times H_{n0}$.

Similarly we can define $m(E_1\times E_2\times\cdots\times E_n)$
as in Theorem \ref{ffcr1} and get the theorem for multi-parameter
integral.
\begin{thm}
Let $F(x_1,\cdots, x_n)$ be a continuous function of bounded
$\Phi_i$-variation in $x_i$ uniformly in
$(x_1,\cdots,x_{i-1},x_{i+1},\cdots, x_n)$; $i=1,\cdots n$,
$G(x_1,\cdots,x_n)$ be of bounded $\Psi_1,\cdots,\Psi_n$-variation
in $(x_1,\cdots,x_n)$ and satisfy the finite large jump condition,
where $\Psi_1,\cdots,\Psi_n$ are convex
functions. If there exist
 monotone increasing functions $\varrho_1,
\cdots, \varrho_n$ subject to $\varrho_1(u)\cdots\varrho_n(u)=u$
such that
\begin{eqnarray}
\sum_{k_n=1}^{\infty}\cdots\sum_{k_1=1}^{\infty}\varrho_1[\varphi_1({1\over
k_1})]\cdots\varrho_n[\varphi_n({1\over k_n})]\psi_1\big[{1\over
{k_1}}[\cdots\psi_n({1\over {k_n}})\cdots]\big]<\infty,
\end{eqnarray}
\end{thm}
 then the integral
\begin{eqnarray*}
&&\int_{a_n}^{b_n} \cdots\int_{a_1}^{b_1}F(x_1,\cdots,x_n)d_{x_1,\cdots,x_n}G(x_1,\cdots,x_n)\nonumber\\
&&\hskip -0.2cm=\lim_{m(E_1\times\cdots E_n)\to 0}
\sum_{k_n=1}^{N_n}\cdots\sum_{k_1=1}^{N_1} F(x_1^{k_1-1},\cdots,
x_n^{k_n-1})(\Delta_{x_n^{k_n-1},x_n^{k_n}}\cdots\Delta_{x_1^{k_1-1},x_1^{k_1}}G)
\end{eqnarray*}
is well defined, as long as $E_1\times E_2\times\cdots\times E_n$
include $H_1\times H_2\times\cdots\times H_n$.

\section{Two-parameter integrals of local times}
\setcounter{equation}{0}

Assume that $X=(X_t)_{t\geq0}$ is a continuous semi-martingale,
$L_t^x$ is the local time of $X_t$ at $x$. By localization
argument, we can assume that $X_t$ is bounded and its local time
$L_t(x)$ is also bounded uniformly in $x$ (see \cite{feng}). We
divide $[0,t]\times [-N,N]$ by $E\times E':=E_{[0,t]}\times
E'_{[-N,N]}=\{0=s_0<s_1<\cdots<s_m=t, -N=x_0<x_1<\cdots<x_l=N\}$,
 where $[-N,N]$ covers the support of local time $L_t^x$.

In this section we will define $\int_{-\infty}^\infty\int_0^t g(s,x) d_{s,x}L_s^x$. We will first use
Theorem \ref{ffcr1} to define the integral
$\int_{-\infty}^\infty\int_0^t \tilde L_s^x d_{s,x}g(s,x)$. Here $\tilde L_s^x$ refers to the continuous part in decompotition (\ref{hz102}) of local times.
\begin{thm}\label{our4}
Assume $g: [0,t]\times R\to R$ is of bounded $\Phi_1,
\Psi_1$-variation in $(s,x)$, i.e. $\sup\limits _{E\times E'}
\sum\limits_{i=0}^{l-1}\Psi_1\big(\sum\limits_{j=0}^{m-1}\Phi_1(|\Delta_j\Delta_i
g|)\big)< \infty$ for the partition we defined as before
and satisfy the finite large jump condition. Then if
there exist monotone increasing functions $\varrho$ and $\sigma$
subject to $\varrho(u)\sigma(u)=u$ such that for a $\delta >0$
\begin{eqnarray}\label{cond4}
\sum_{n,m}\varrho[({1\over n})^{1\over {2+\delta}}]\sigma({1\over
m})\varphi_1({1\over n}\psi_1({1\over m}))<\infty,
\end{eqnarray}
the integral
\begin{eqnarray}
&&\int_{-\infty}^{\infty}\int_0^t \tilde L_s^x d_{s,x}g(s,x)\nonumber\\
&=&\lim_{m(E\times E')\to0}\sum _{i=0}^{l-1}\sum_{j=0}^{m-1}\tilde L(s_j,x_i)\Big(g(s_{j+1}, x_{i+1})-g(s_{j+1}, x_{i})\nonumber\\
&&\hskip 3cm-g(s_{j}, x_{i+1})+g(s_{j}, x_{i})\Big)
\end{eqnarray}
is well defined for almost all $\omega\in\Omega$ in the sense of
Theorem \ref{ffcr1}.
\end{thm}
{\bf Proof:} Note ${L}_s(x)$ is increasing in $s$ so  of
bounded variation in $s$. Let $h$ be defined by (\ref{hz103}). It
is easy to know from (\ref{hz15}), $h(s,x)$ is of bounded
variation in $s$. Therefore we have $\tilde L(s,x)$ is of bounded
variation in $s$. In particular, using (\ref{hz102}),
(\ref{hz103}) and (\ref{hz15}) we obtain
\begin{eqnarray*}
\sup _{E} \sum_{j=0}^{l-1}| \tilde L(s_{j+1}, x)-\tilde L(s_{j}, x) |\leq L_t(x)+\int _0^t|dV_s|\leq P,
\end{eqnarray*}
where $P$ is a constant independent of x. On the other hand, from
Lemma \ref{newlem} and Lemma \ref{newlem1}, we know,
\begin{eqnarray}
\sup _{E'}\sum_{i=0}^{m-1}|\tilde L(s, x_{i+1})-\tilde L(s,
x_{i})|^{2+\delta}\nonumber \leq Q,
\end{eqnarray}
where $Q$ is a constant independent of $s$. Therefore under
condition (\ref{cond4}), the integral $\int_{-\infty}^\infty
\int_0^t \tilde L_s(x)d_{s,x}g(s,x)$ is well defined. $\hfill\diamond$
\bigskip

\begin{cor}\label{our5}
Assume $g: [0,t]\times R\to R$ is of bounded $p,q$-variation, i.e. $\sup\limits_{E\times E'}\sum\limits_{i=0}^{l-1}\left(\sum\limits_{j=0}^{m-1}|\Delta_j \Delta_i g|^p\right)^q< \infty$, where $p,q \geq 1$, $2q+1>2pq$ and satisfies the finite large jump condition, then the integral
\begin{eqnarray}\label{local}
&&\int_{-\infty}^{\infty}\int_0^t \tilde L_s^x d_{s,x}g(s,x)\nonumber\\
&=&\lim_{m(E\times E')\to0}\sum _{i,j}\tilde L(s_j,x_i)\Big(g(s_{j+1}, x_{i+1})-g(s_{j+1}, x_{i})\nonumber\\
&&\hskip 3.5cm-g(s_{j}, x_{i+1})+g(s_{j}, x_{i})\Big)
\end{eqnarray}
is well defined in the sense of Theorem \ref{ffcr1}.
\end{cor}
{\bf Proof:} For any $p,q \geq 1$ satisfying $2q+1>2pq$, we have
$2(1-{1\over p})<{1\over {pq}}$. Therefore there exists a number
$\alpha$ such that $2(1-{1\over p})<\alpha<{1\over {pq}}$. This
implies that ${\alpha\over 2}+{1\over p}>1$ and ${1-\alpha}
+{1\over pq}>1$. So there is a small $\delta>0$ such that
${\alpha\over {2+\delta}}+{1\over p}>1$ and ${1-\alpha} +{1\over
pq}>1$. Take $\varrho(u)=u^{\alpha}$ and $\sigma(u)=u^{1-\alpha}$,
then it is easy to see that
\begin{eqnarray}\label{b1}
\sum_{n,m}\varrho[({1\over n})^{1\over {2+\delta}}]\sigma ({1\over
m})({1\over n})^{1\over p} ({1\over m})^{1\over pq}<\infty.
\end{eqnarray}
Therefore the integral (\ref{local}) is well defined.
$\hfill\diamond$
\vskip5pt

After defining the integral
$\int_{-\infty}^\infty \int_0^t\tilde L_s(x)d_{s,x}g(s,x)$, let's study
the integral $\int_{-\infty}^\infty \int_0^tg(s,x)d_{s,x}\tilde L_s^x$.
Note
\begin{eqnarray}\label{e1}
&&\sum_{i=0}^{l-1}\sum_{j=0}^{m-1} g(s_j,x_i)\left[\tilde L_{s_{j+1}}(x_{i+1})-\tilde L_{s_{j}}(x_{i+1})-\tilde L_{s_{j+1}}(x_{i})+\tilde L_{s_{j}}(x_{i})\right]\nonumber\\
&=& \sum_{i=1}^{l}\sum_{j=1}^{m} g(s_{j-1},x_{i-1})\tilde L_{s_j}(x_i)-\sum_{i=1}^{l}\sum_{j=0}^{m-1} g(s_{j},x_{i-1})\tilde L_{s_j}(x_i)\nonumber\\
&&-\sum_{i=0}^{l-1}\sum_{j=1}^{m} g(s_{j-1},x_{i})\tilde L_{s_j}(x_i)+\sum_{i=0}^{l-1}\sum_{j=0}^{m-1} g(s_{j},x_{i})\tilde L_{s_j}(x_i)\nonumber\\
&=& \sum_{i=1}^{l}\sum_{j=1}^{m} \tilde L_{s_j}(x_i)\left[g(s_j,x_i)-g(s_j,x_{i-1})-g(s_{j-1},x_i)+g(s_{j-1},x_{i-1})\right]\nonumber\\
&&-\sum_{i=1}^l \left[g(0,x_{i-1})\tilde L_0(x_i)-g(t,x_{i-1})\tilde L_t(x_{i})\right]\nonumber\\
&&-\sum_{j=1}^m\left[g(s_{j-1},-N)\tilde L_{s_j}(-N)-g(s_{j-1},N)\tilde L_{s_j}(N)\right]\nonumber\\
&&+\sum_{j=0}^{m-1}\left[g(s_j,-N)\tilde L_{s_j}(x_0)-g(s_j,N)\tilde L_{s_j}(N)\right]\nonumber\\
&&+\sum_{i=0}^{l-1}\left[g(0,x_i)\tilde L_0(x_i)-g(t,x_i)\tilde L_t(x_i)\right]\nonumber\\
&=& \sum_{i=1}^{l}\sum_{j=1}^{m} \tilde L_{s_j}(x_i)\left[g(s_j,x_i)-g(s_j,x_{i-1})-g(s_{j-1},x_i)+g(s_{j-1},x_{i-1})\right]\nonumber\\
&&-\sum_{i=1}^l\tilde L_t(x_i)(g(t,x_i)-g(t,x_{i-1})).
\end{eqnarray}
Under the conditions of Theorem \ref{our4} and Proposition
\ref{proposition} and noticing that $\tilde L_t(x)$ is continuous
in $t$, we know that the first term of (\ref{e1}) converges to
$\int_{-\infty}^\infty\int_0^t \tilde L_s(x)d_{s,x}g(s,x)$, and
from rough path integration of one parameter, we know that the
second term converges to $\int_{-\infty}^\infty \tilde
L_t(x)d_xg(t,x)$ if further $g(s,x)$ is of bounded
$\gamma-$variation ($1\leq \gamma <2$) in $x$ uniformly in $s$. So
the sum
\begin{eqnarray*}
\sum\limits_{i=0}^{l-1}\sum\limits_{j=0}^{m-1}
g(s_j,x_i)\left[\tilde L_{s_{j+1}}(x_{i+1})-\tilde L_{s_{j}}(x_{i+1})-\tilde L_{s_{j+1}}(x_{i})+\tilde L_{s_{j}}(x_{i})\right]
\end{eqnarray*}
converges, we denote its limit by
\begin{eqnarray}\label{new}
&&\int_{-\infty}^\infty\int_0^t g(s,x)d_{s,x}\tilde L_s^x\nonumber\\
&=&\lim_{m(E\times E')\to 0}\sum_{i=0}^{l-1}\sum_{j=0}^{m-1} g(s_j,x_i)\Big[\tilde L_{s_{j+1}}(x_{i+1})-\tilde L_{s_{j}}(x_{i+1})\nonumber\\
&&\hskip 3cm-\tilde L_{s_{j+1}}(x_{i})+\tilde L_{s_{j}}(x_{i})\Big],
\end{eqnarray}
and
\begin{eqnarray}
\hskip -0.5cm\int_{-\infty}^\infty\int_0^t
g(s,x)d_{s,x}\tilde L_s^x=\int_{-\infty}^\infty\int_0^t\tilde L_s^xd_{s,x}g(s,x)-\int_{-\infty}^\infty
\tilde L_t(x)d_xg(t,x).
\end{eqnarray}

Now recall decomposition (\ref{hz102}) and (\ref{hz103}) and Lemma
\ref{newlem1}, as in  Elworthy, Truman and Zhao \cite{Zhao}, the
integral $\int_0^t\int _{-\infty}^{\infty}g(s,x)d_{s,x}h(s,x)$ is
defined as a two-parameter Lebesgue-Stieltjes integral. Therefore
we can define
\begin{eqnarray*}
\int_0^t\int _{-\infty}^{\infty}g(s,x)d_{s,x}L(s,x)=\int_0^t\int
_{-\infty}^{\infty}g(s,x)d_{s,x}\tilde L(s,x) +\int_0^t\int
_{-\infty}^{\infty}g(s,x)d_{s,x}h(s,x).
\end{eqnarray*}
\begin{rmk}\label{fcr7}
If $g(s,x)$ is $C^1$ in $x$, we have
\begin{eqnarray*}
\int_{-\infty}^\infty\int_0^t g(s,x)d_{s,x}L_s^x=-\int_{-\infty}^\infty\int_0^t \nabla g(s,x)d_sL_s(x)dx.
\end{eqnarray*}
This can be seen from the following. As one can always add some
points in the partition to make $L_{s_j}^{x_1}=0$ and
$L_{s_j}^{x_{l+1}}=0$ for all $j=1,2,\cdots,m, $ as $L$ has a
compact support in $x$, therefore
\begin{eqnarray*}
&&\lim_{m(E\times E')\to 0}\sum _{i=1}^{l}\sum _{j=1}^{m} g(s_j,x_i)\Big[L_{s_j}^{s_{j+1}}(x_{i+1})-L_{s_j}^{s_{j+1}}(x_{i})\Big]\\
&=&\lim_{m(E\times E')\to 0}\left(\sum _{i=1}^{l}\sum _{j=1}^{m} g(s_j,x_{i})L_{s_j}^{s_{j+1}}(x_{i+1})-\sum _{i=0}^{l-1}\sum _{j=1}^{m} g(s_j,x_{i+1})L_{s_j}^{s_{j+1}}(x_{i+1})\right)\\
&=&-\lim_{m(E\times E')\to 0}\sum _{i=1}^{l}\sum _{j=1}^{m} \Big[g(s_j,x_{i+1})-g(s_j,x_i)\Big]L_{s_j}^{s_{j+1}}(x_{i+1})\\
&=&-\lim_{m(E\times E')\to 0}\sum _{i=1}^{l}\sum _{j=1}^{m}\nabla g(s_j,\xi_i)L_{s_j}^{s_{j+1}}(x_{i+1})(x_{i+1}-x_i)\\
&=&-\lim_{m(E'_{[-N,N]})\to 0}\sum_i\int_0^t \nabla g(s,x_{i+1})d_sL_s(x_{i+1})(x_{i+1}-x_i)\\
&&-\lim_{m(E'_{[-N,N]})\to 0}\sum_i\int_0^t \bigg(\nabla g(s,\xi_i)-\nabla g(s,x_{i+1})\bigg) d_sL_s(x_{i+1})(x_{i+1}-x_i)\\
&=&-\int_{-\infty}^\infty\int_0^t \nabla g(s,x)d_sL_s(x)dx.
\end{eqnarray*}

\end{rmk}

\begin{thm}\label{fcr5}
Let $f: [0,t]\times R\to R$ be of bounded $\gamma$-variation in $x$ uniformly in $s$ and of bounded $p, q$-variation in $(s,x)$ and satisfy the finite large jump condition, where $1\leq \gamma <2$ and $p,q \geq 1$, $2q+1>2pq$, and
\begin{eqnarray}
f_n(s,x):=\int _0^2\int _0^2\rho(r)\rho(z)f(s-{r\over
n},x-{z\over n})drdz, \ \ n\geq 1,
\end{eqnarray}
where $\rho$ is the mollifier defined in (\ref{mollifier}).
 Then
\begin{eqnarray*}
\int_{-\infty}^{\infty}\int_0^t f_n(s,x)d_{s,x}L_s^x\to\int_{-\infty}^{\infty}\int_0^t f(s,x)d_{s,x}L_s^x,\ as\ n\to \infty.
\end{eqnarray*}
\end{thm}
{\bf Proof}: First we can easily verify that $f_n$ are also of bounded $p,q$-variation.
 We extend $f$ to $s<0$ by defining $f(s,x)\equiv 0$, for $s<0$, and denote an arbitrary partition of $[0,t]\times [-N-2, N]$ by
\begin{eqnarray*}
E\times E'_1:=\{0=s_0<s_1<\cdots<s_{m}=t, -N-2=x_0<x_1<\cdots<x_{l'}=N\}.
\end{eqnarray*}
Note $[-N-2,N]$ also covers the compact support of local time, and
\begin{eqnarray*}
\sup\limits_{E\times
E'_1}\sum_{i=1}^{l'}\left(\sum_{j=1}^{m}|\Delta_j\Delta_if|^p\right)^q=M,
\end{eqnarray*}
and
\begin{eqnarray*}
\sup\limits_{E'_1}\sum_{i=1}^{l'}
|f(s,x_i)-f(s,x_{i-1})|^\gamma=M',
\end{eqnarray*}
where $M$ and $M'$ are constants. So by H\"older inequality,
\begin{eqnarray*}
&&\sum_{i=1}^l\left(\sum_{j=1}^m|\Delta_j\Delta_if_n|^p\right)^q\\
&=&\sum_{i=1}^l \left(\sum_{j=1}^m \left|\int_0^2 \int_0^2\rho(r)\rho(z)\Delta_j\Delta_if(\cdot-{r\over n}, \cdot-{z\over n})drdz\right|^p\right)^q\\
&\leq&A\sum_{i=1}^l\left(\int_0^2\int_0^2\sum_{j=1}^m\left|\Delta_j\Delta_i f(\cdot-{r\over n},\cdot-{z\over n})\right|^p drdz\right)^q\\
&\leq&B\int_0^2\int_0^2\sum_{i=1}^l\left(\sum_{j=1}^m\left|\Delta_j\Delta_i f(\cdot-{r\over n},\cdot-{z\over n})\right|^p\right)^qdzdr\\
&\leq &B\int_0^2\int_0^2\sup\limits_{E\times E'_1}\sum_{i=1}^{l'}\left(\sum_{j=1}^{m}|\Delta_j\Delta_if|^p\right)^qdrdz\\
&\leq & M_1,
\end{eqnarray*}
where $A$, $B$ and $M_1$ (independent of $n$) are constants.
Also from the above estimate, the finite large jump condition for $f_n$ when $n$ is sufficiently large
follows from the finite large jump assumption of $f$.
Similarly,
\begin{eqnarray*}
&&\sum_{i=1}^l|f_n(s, x_i)-f_n(s,x_{i-1})|^\gamma\\
&=&\sum_{i=1}^l|\int_0^2\int_0^2\rho(r)\rho(z)\left(f(s-{r\over n}, x_i-{z\over n})-f(s-{r\over n}, x_{i-1}-{z\over n}\right)drdz|^\gamma\\
&\leq& C \int_0^2 \int_0^2\sum_{i=1}^l|f(s-{r\over n}, x_i-{z\over n})-f(s-{r\over n}, x_{i-1}-{z\over n})|^\gamma drdz\\
&\leq& C \int_0^2 \int_0^2\sup\limits_{E'}\sum_{i=1}^{l'} |f(s,x_i)-f(s,x_{i-1})|^\gamma drdz\\
&\leq& M_2
\end{eqnarray*}
where $C$ and $M_2$ (independent of $n$) are constants. So the
integral\\ $\int_{-\infty}^\infty\int_0^t f_n(s,x)d_{s,x}L_s^x$ is
well defined, by argument we discussed before,
\begin{eqnarray}\label{hz104}
&&
\int_{-\infty}^\infty\int_0^t
f_n(s,x)d_{s,x}L_s^x\nonumber\\
&=&\int_{-\infty}^\infty\int_0^t
\tilde L_s^xd_{s,x}f_n(s,x)-\int_{-\infty}^\infty \tilde L_t^xd_xf_n(t,x)\nonumber\\
&&
+\int_{-\infty}^\infty\int_0^t
f_n(s,x)d_{s,x}h(s,x).
\end{eqnarray}
 For
such $p, q$ satisfying $p,q\geq 1$, and $2q(p-1)<1$, there exist a
small positive number $\delta >0$ such that
$(2+\delta)q(p+\delta-1)<1$, so
\begin{eqnarray*}
\sum_{n,m}\varrho[({1\over n})^{1\over {2+\delta}}]\sigma ({1\over
m})({1\over n})^{1\over {p+\delta}} ({1\over m})^{1\over
({p+\delta})q}<\infty
\end{eqnarray*}
still holds for $\rho(u)=u^\alpha$, $\sigma(u)=u^{1-\alpha}$,
where $(2+\delta)(1-{1\over {p+\delta}})<\alpha<{1\over
{(p+\delta)q}}$. By Theorem \ref{our3} and Proposition
\ref{proposition}, we can pass the limit to get
\begin{eqnarray*}
\lim_{n\to\infty}\int_{-\infty}^\infty\int_0^t \tilde L_s^xd_{s,x}f_n(s,x)=\int_{-\infty}^\infty\int_0^t \tilde L_s^xd_{s,x}f(s,x).
\end{eqnarray*}
Using a similar method as in the proof of Theorem \ref{fcr1}, we can prove that
\begin{eqnarray*}
\lim_{n\to\infty}\int_{-\infty}^\infty \tilde L_t^xd_xf_n(t,x)=\int_{-\infty}^\infty \tilde L_t^xd_xf(t,x).
\end{eqnarray*}
The convergence of the  last term $\int_{-\infty}^\infty\int_0^t
f_n(s,x)d_{s,x}h(s,x)$ in (\ref{hz104}) follows from Lebesgue's dominated convergence
theorem. So
we proved the desired result.
$\hfill\diamond$ \vskip5pt

 \begin{thm}\label{last}
 Let $X=(X_s)_{s\geq0}$ be a continuous semi-martingale and assume
 $f: [0,\infty)\times R\to R$ satisfies
\vskip2pt

(i) $f$ is left continuous and locally bounded, with $f(t,x)$ jointly continuous from the right in $t$
and left in $x$ at each point $(0,x)$,
\vskip2pt

(ii) the left derivatives ${\partial ^-\over \partial t}f$ and  $\nabla ^-f$ exist at all points of $(0,
\infty )\times R$ and $[0,\infty)\times R$ respectively,
\vskip2pt

(iii)  ${\partial ^-\over \partial t}f$ and  $\nabla ^-f$ are left continuous and locally bounded,
\vskip2pt

(iv)  $\nabla ^-f(t,x)$ is of bounded $\gamma$-variation in $x$ uniformly in $t$ and of bounded $p, q$-variation in ($t,x$) and satisfies the finite large jump condition, where $1\leq \gamma<2$, and $p,q \geq 1$, $2q+1>2pq$.
\vskip2pt

\noindent
Then we have:
 \begin{eqnarray}
 f(t,X_t)&=&f(0, X_0)+\int_0^t{\partial^-\over{\partial s}}f(s,X_s)ds+\int_0^t \nabla ^-f(s,X_s)dX_s\nonumber\\
 &&-\int_0^t\int_{-\infty}^{\infty}\nabla^-f(s,x)d_{s,x} L_s^x,
 \end{eqnarray}
 where $L_t^x$ is the local time of $X_t$ at $x$, the last integral is defined in (\ref{new}).
  \end{thm}
{\bf Proof}: Similar to the proof in \cite{Zhao}, we can use
smoothing procedure and take the limit to prove our result. The
main different key point is the following : by Remark \ref{fcr7}
and Theorem \ref{fcr5},
\begin{eqnarray*}
 &&{1\over2}\int_0^t \Delta f_n(s,X_s)d<M>_s\\
&=&\int_{-\infty}^\infty\int_0^t \Delta f_n(s,x) dL_s^x dx\\
&=&-\int_{-\infty}^\infty\int_0^t \nabla f_n(s,x) d_{s,x}L_s^x\\
&\to& -\int_{-\infty}^\infty \int_0^t\nabla^- f(s,x) d_{s,x}L_s^x,
 \end{eqnarray*}
 when $n\to \infty$.
 $\hfill\diamond$
 \vskip5pt

 \begin{exap} Consider a function $f(t,x)=x^3t^3\cos ({1\over t}+{1\over x})$ for $t,x\neq 0$ and $f(t,0)=f(0,x)=f(0,0)=0$.
This function is $C^{1,1}$ and its derivative about $x$ is
${\partial\over
\partial x}f(t,x)=3t^3x^2\cos ({1\over t}+{1\over x})+xt^3\sin ({1\over t}+{1\over x})$
for $t,x\neq 0$ and ${\partial\over
\partial x}f(t,0)={\partial\over
\partial x}f(0,x)={\partial\over
\partial x}f(0,0)=0$. It is easy to
see that ${\partial\over
\partial x}f(t,x)$ is of unbounded variation in $x$ and in $(t,x)$ , but
of $\gamma$-variation in $x$ for any $\gamma>1$, $p,1$-variation
in $(t,x)$ for any $p>1$ (similar to Example 3.1). So Theorem
\ref{last} can be used.
\end{exap}

 Finally we would like to mention that our result should also work
 for stable processes noticing the $p$-variation result on the local times of stable
 processes studied by Marcus and Rosen \cite{rosen}. But we should also point out
 that Marcus and Rosen's definition to $p$-variation is different from ours.
 But we can extend the proof of Lemma \ref{newlem} to stable
 process. For the length of the paper, these results are not
 included in this paper.

\bigskip

{\bf Acknowledgement} \vskip5pt 

We would like to thank N.
Eisenbaum, K.D. Elworthy, T. Lyons, Z.M. Ma, Z. Qian, T. Zhang and S. Peng for useful conversation,
especially to T. Lyons who introduced the idea of Young integral
to us during the UK-Japan Winter School in January 2003 at Warwick
University and to S. Peng for his hospitality and useful
conversations during HZ's visit to Shandong University where the
paper was finalized. EPSRC's grant no. GR/R69518 is gratefully
acknowledged.

\end{document}